\documentclass{article}
\usepackage{graphicx}
\usepackage{pgfplots}
\usepackage{subfigure}
\usepackage{caption}
\usepackage[justification=centering]{caption}
\usepackage{amsmath}
\usepackage{amssymb}
\usepackage{amsthm}
\usepackage{enumerate}
\usepackage{listings}
\usepackage{epstopdf}
\usepackage{float}
\usepackage[utf8]{inputenc}
\usepackage{enumitem}
\usepackage[top=2cm, bottom=2cm, left=2cm, right=2cm]{geometry}
\newtheorem{theorem}{Theorem}
\newtheorem{lemma}{Lemma}
\usepackage{array}
\usepackage{authblk}

\usepackage{cite}

\usepackage{setspace}
\begin{document}
	
\title{\textbf{Dynamics Motivated by Sierpinski Fractals}}
\author[1]{\textbf{Marat Akhmet}\thanks{Corresponding Author Tel.: +90 312 210 5355, Fax: +90 312 210 2972, E-mail: marat@metu.edu.tr}$^{,}$}
\author[2]{\textbf{Mehmet Onur Fen}}
\author[1]{\textbf{Ejaily Milad Alejaily}}
\affil[1]{\textbf{Department of Mathematics, Middle East Technical University, 06800 Ankara, Turkey}}
\affil[2]{\textbf{Department of Mathematics, TED University, 06420 Ankara, Turkey}}

\date{}
\maketitle
	
\begin{center}
	\textbf{\large Abstract}
\end{center}
Fatou-Julia iteration (FJI) is an effective instrument to construct fractals. Famous Julia and Mandelbrot sets are strong confirmations of this. In the present study, we use the paradigm of FJI to construct and map Sierpinski fractals. The fractals can be mapped by developing a mapping iteration on the basis of FJI. Because of the close link between mappings, differential equations and dynamical systems, one can introduce dynamics for a fractal through differential equations such that it becomes points of the solution trajectory. Thus, in this paper, we consider two types of dynamics motivated by the Sierpinski fractals. The first one is the dynamics of FJI itself, and the second one is the dynamics of a fractal mapping iteration which can be performed through differential equations. The characterization of fractals as trajectory points of the dynamics can help to enhance and widen the scope of their applications in physics and engineering. \\


\section{\normalfont INTRODUCTION}
The term ``\textit{fractal}" was coined by Benoit Mandelbrot in 1975 \cite{Mandelbrot0}. He defined fractal as a set for which the Hausdorff dimension strictly exceeds the topological dimension \cite{Mandelbrot1}. Self-similarity and fractional dimension are the most two important features of fractals. The connection between them is that self-similarity is the easiest way to construct a set that has fractional dimension \cite{Crownover}.

Dealing with fractals goes back to the 17th Century when Gottfried Leibniz introduced the notions of recursive self-similarity \cite{Zmeskal}. A considerable leap in the construction of fractals was performed in 1883 by Georg Cantor as he discovered the most essential and influential fractal known as the Cantor set.

\begin{figure}[]
	\centering
	\subfigure[Sierpinski gasket ]{\includegraphics[width = 2.5in]{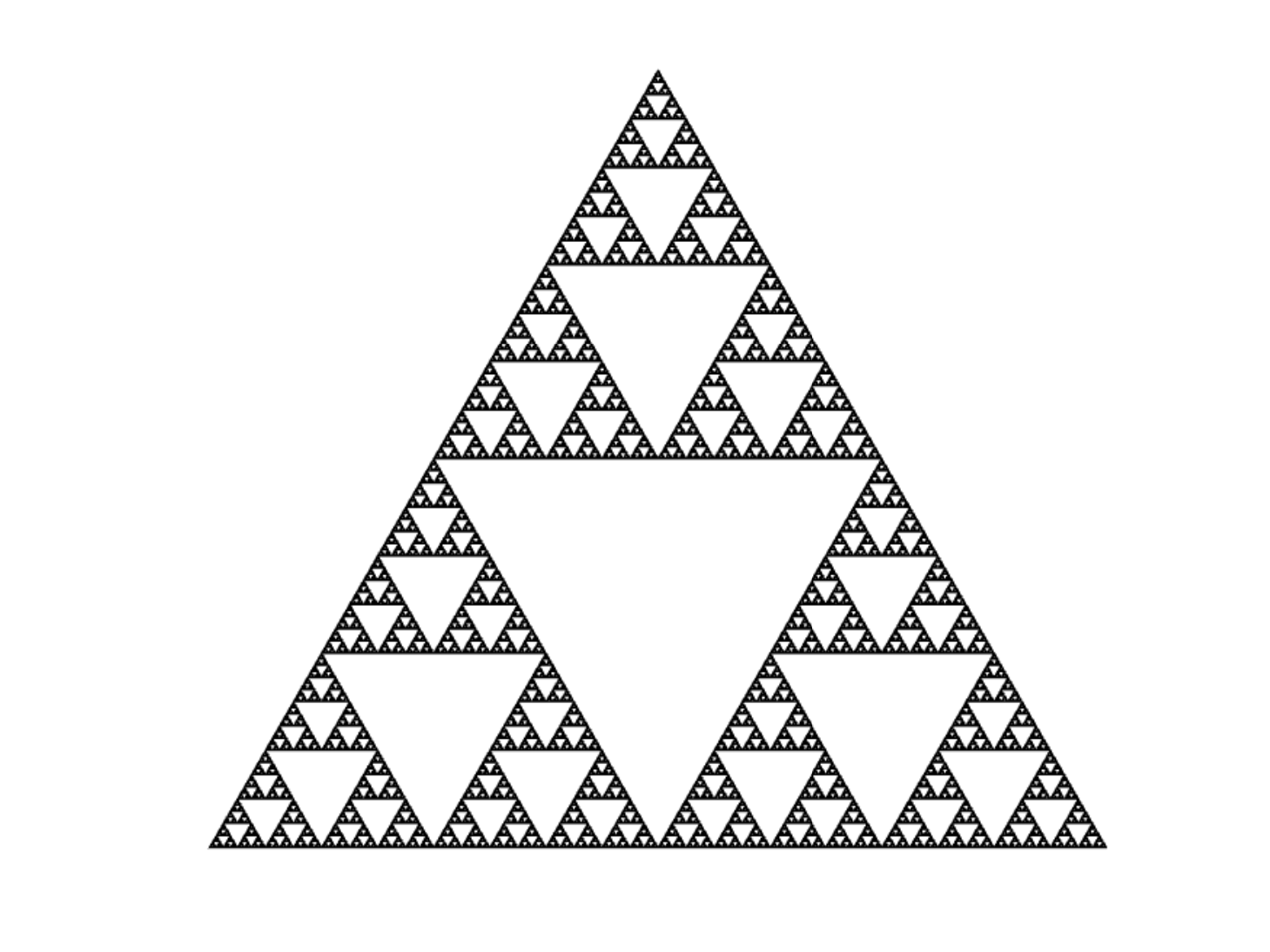}\label{SFG}} 
	\subfigure[Sierpinski carpet]{\includegraphics[width = 2.5in]{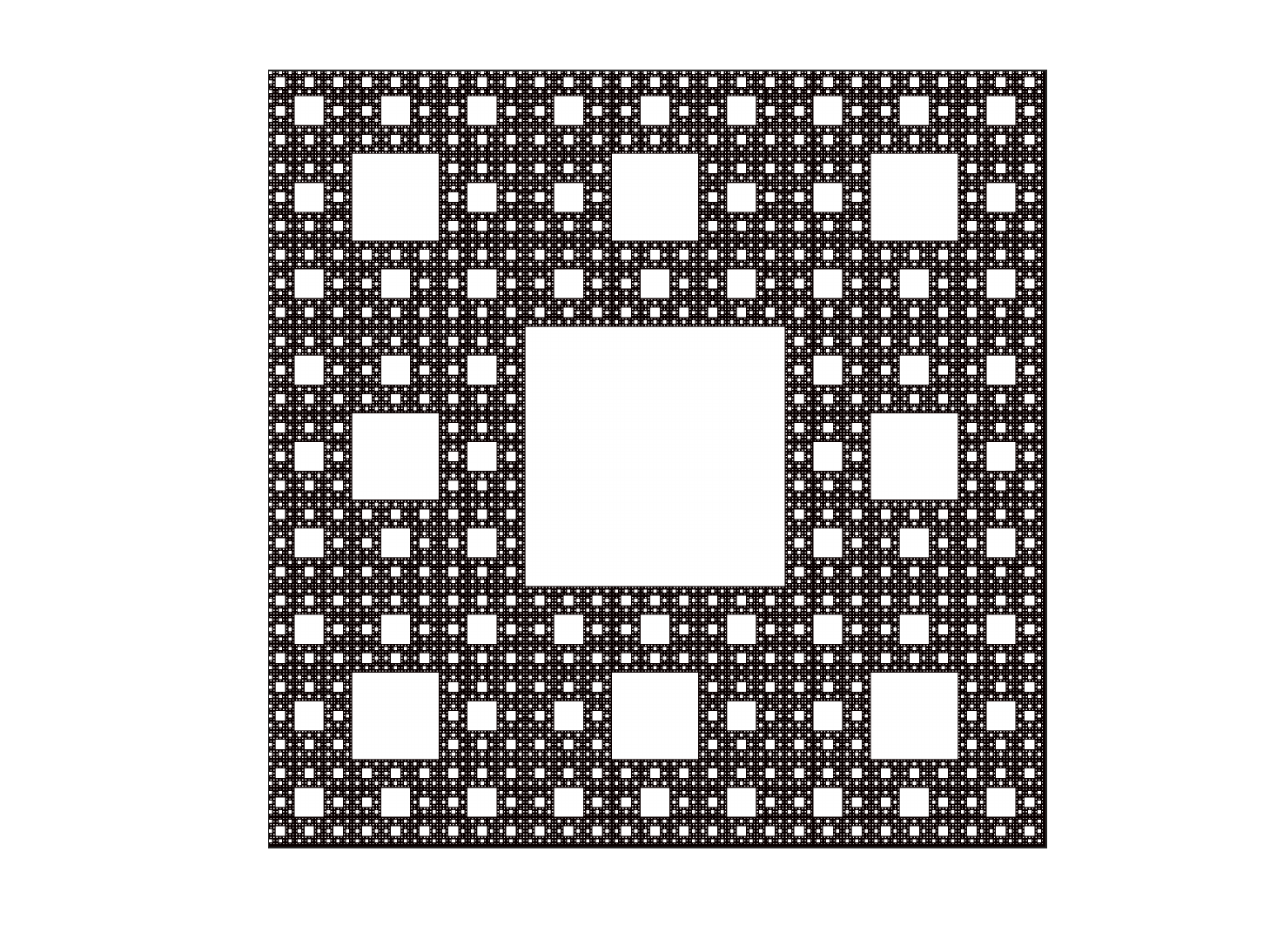}\label{SFC}}
	\caption{Sierpinski Fractals}
	\label{SierpinskiFractals}   				
\end{figure}

Waclaw Sierpinski was one of the mathematicians who made significant contributions in the field of fractals. He introduced the famous triangular fractal in 1916, known as the Sierpinski gasket. The fractal is generated by a recursive process of removing symmetrical parts from an initial triangle. The process starts with subdividing an equilateral solid triangle into four identical sub-triangles then removing the central one. In the next iterations, the same procedure is repeated to each of the remaining triangles from the preceding iteration. In an analogous way to the gasket, Sierpinski developed a square fractal known as the Sierpinski carpet. Figure \ref{SierpinskiFractals} shows the two Sierpinski fractals for finite iterations, which are constructed by the methods introduced and discussed in the present paper.

Involvement of the dynamics of iterative maps in fractal construction was a critical step made by the French mathematicians Pierre Fatou and Gaston Julia around 1917-1918, during their independent studies on the iteration of rational functions in the complex plane \cite{Julia,Fatou}. They described what we call today the Fatou-Julia iteration (FJI) \cite{Mandelbrot2}. The iteration is defined over a domain $ \mathcal{D} \subseteq \mathbb{C} $ by
\begin{equation} \label{FJIt}
z_{n+1}=F(z_n),
\end{equation}
where $ F: \mathcal{D} \to \mathcal{D} $ is a given function for the construction of the fractal set $ \mathcal{F} $. The points $ z_0 \in \mathcal{D} $ are included in the set $ \mathcal{F} $ depending on the boundedness of the sequence $ \{z_n\}, \; n=0, 1, 2, ... $, and we say that the set $ \mathcal{F} $ is constructed by FJI.

In practice one cannot verify the boundedness for infinitely long iterations. This is why in simulation we fix an integer $ k $ and a bounded subset $ M \subset \mathbb{C} $, and denote by $ \mathcal{F}_k $ the collection of all points $ z_0 \in \mathcal{D} $ such that the points $ z_n $ where the index $ n $ is between $ 1 $ and $ k $, $ n=1, 2, ... , k $, belong to $ M $. In what follows we call the set $ \mathcal{F}_k $ the \textit{kth approximation} of the set $ \mathcal{F} $.

The most popular fractals, Julia and Mandelbrot sets, are generated using the iteration of the quadratic map  $ F(z_n)=z_n^2+c $, where $ c $ is a complex number. Julia sets contain the points $ z_0 \in \mathbb{C} $  corresponding to the bounded sequence $ \{z_n\} $, whereas the Mandelbrot set is the set of parameter values  $ c \in \mathbb{C} $ such that the sequence $ \{z_n\} $, $ z_0=0 $ remains bounded.

In paper \cite{Akhmet} we suggested a new method that makes it possible to map a fractal attained by FJI. The method is based on involving the map $ f $ in the iteration (\ref{FJIt}). The Fractal Mapping Iteration (FMI) is defined as
\begin{equation} \label{FMIt}
f^{-1}(z_{n+1})=F \big(f^{-1}(z_n) \big).
\end{equation}
In this recursive equation we apply FJI to the preimage of the mapped set. Therefore, the condition of boundedness should be applied to the sequence $ \{ f^{-1}(z_n) \} $. The FMI (\ref{FMIt}) is described in general form, i.e. it is valid for any function $ F $. In paper \cite{Akhmet} we applied the FMI to specific fractals, namely Julia and Mandelbrot sets, and constructed discrete and continuous dynamics for them using differential equations.

In the present paper we develop the application of FMI for Sierpinski fractal sets. For that purpose, in the upcoming section we introduce some new schemes for the construction of Sierpinski carpet and gasket. Sections 3 and 4 are devoted to mapping of Sierpinski fractals and constructing continuous dynamics for these sets utilizing FMI and differential equations.
	
\section{\normalfont CONSTRUCTION OF FRACTALS}

The Sierpinski fractals are typically generated using iterated function system \cite{Hutchinson,Barnsley}, which is defined as a collection of affine transformations. Fatou-Julia iteration (sometimes called ``Escape Time Algorithm") can be constructed from IFS \cite{Barnsley,Barnsley1}. In this section, we adopt the idea of FJI and develop some schemes for constructing Sierpinski fractals. The technique of the FJI is based on detecting the points of a fractal set through the boundedness of their iterations under a specific map. Here, we shall extend the technique to include any possible criterion for grouping points in a given domain.

\subsection{\normalsize Sierpinski carpets}	

At first glance, the Sierpinski carpet seems to be a two dimensional version of the middle third Cantor set. To discuss this thought, let us consider the tent map $ T $ defined on the interval $ \mathcal{I}=[0, 1] $ such that
\begin{equation} \label{TentMap}
T(x) = \left\{ \begin{array}{ll}\vspace{2mm}
3 \, x,  & \; \text{if} \;\; x \leq \frac{1}{2}, \\ 
3 (1-x), & \; \text{if} \;\; x > \frac{1}{2}.
\end{array}
\right.
\end{equation}
The Julia set corresponding to the map (\ref{TentMap}) is the middle third Cantor set \cite{Mandelbrot2}. For planar fractals we  consider the FJI defined by the two-dimensional tent map 
\begin{equation} \label{2DTentMap}
(x_{n+1}, y_{n+1}) = \bigg(\frac{3}{2}-3 \big|x_n-\frac{1}{2}\big|, \; \frac{3}{2}-3 \big|y_n-\frac{1}{2}\big|\bigg),
\end{equation}
with the initial square $ \mathcal{D}=[0, 1]\times[0, 1] $. If we exclude each point $ (x_0, y_0) $  whose iterated values $ (x_n, y_n) $ escape from $ \mathcal{D} $, i.e., at least one coordinate escapes, $ x_n > 1 $ or $ y_n > 1 $ for some $ n \in \mathbb{N} $, we shall get a Cantor dust. This set is simply the Cartesian product of the Cantor set with itself, and it is a fractal possessing both self-similarity and fractional dimension. Figure \ref{CD,pScs} (a) shows the $ 3 $rd approximation of the Cantor dust generated by (\ref{2DTentMap}).

Let us now modify this procedure such that a point $ (x_0, y_0) $ is excluded from $ \mathcal{D} $ if both of its coordinates' iterations $ (x_n, y_n) $ escape from the initial set, that is, if  $ x_n > 1 $ and $ y_n > 1 $ for some $ n \in \mathbb{N} $. This procedure for iteration (\ref{2DTentMap}) will give a kind of two dimensional Cantor set shown in Fig. \ref{CD,pScs} (b) with the 3rd approximation. More similar object to the Sierpinski carpet can be obtained by considering simultaneous escape of both coordinates, viz., $ (x_0, y_0) $ is excluded only if  $ x_n > 1 $ and $ y_n > 1 $ at the same iteration $ n $. Figure \ref{CD,pScs} (c) shows the $ 5 $th  approximations of the resulting set. This set is clearly not a fractal from the dimension point of view. The self-similarity is also not satisfied over the whole set. However, a special type of self-similarity can be observed where the corners replicate the whole shape.

\begin{figure}[]
	\centering
	\subfigure[]{\includegraphics[width = 1.8in]{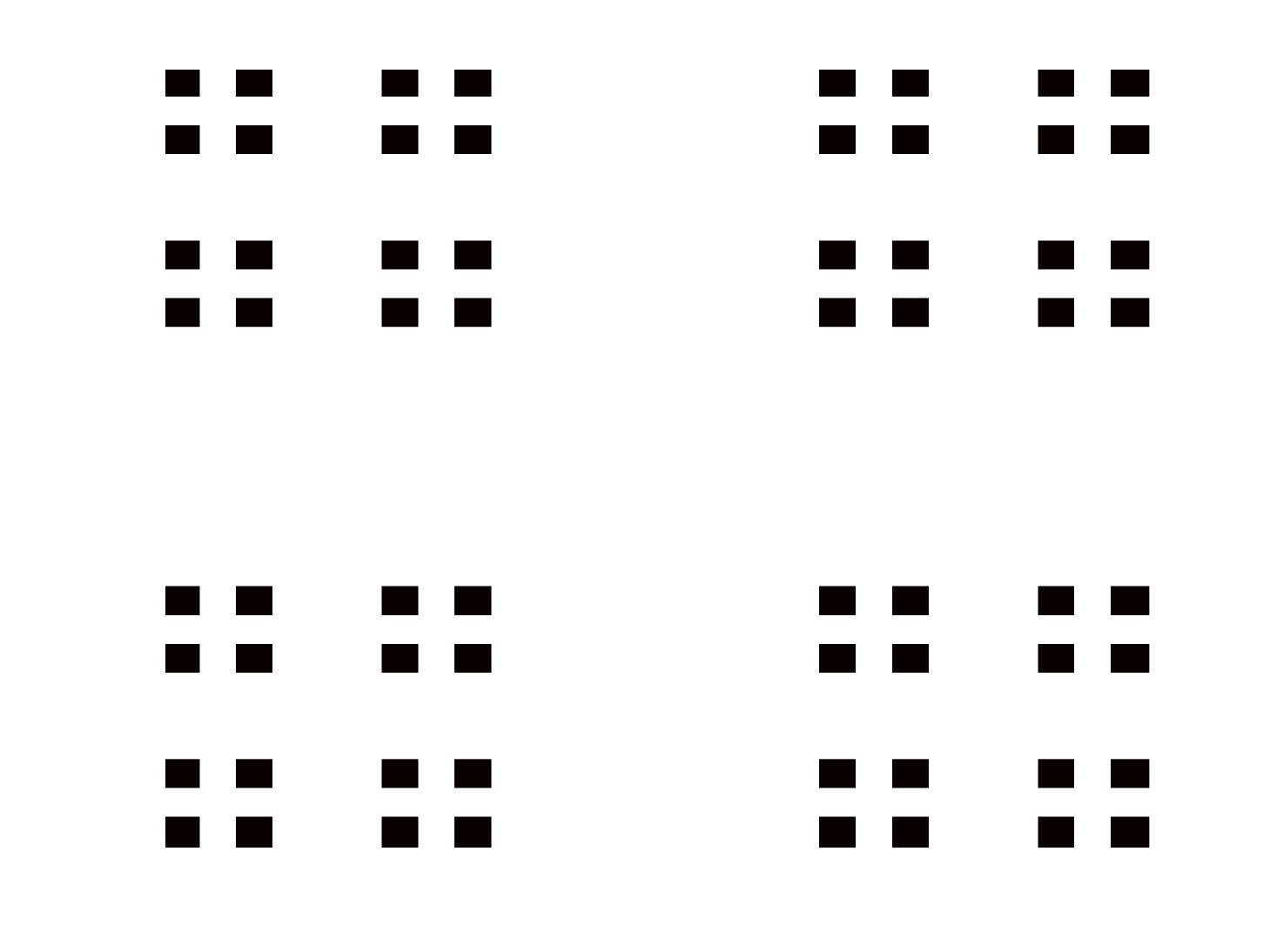}\label{}} 
	\subfigure[]{\includegraphics[width = 1.8in]{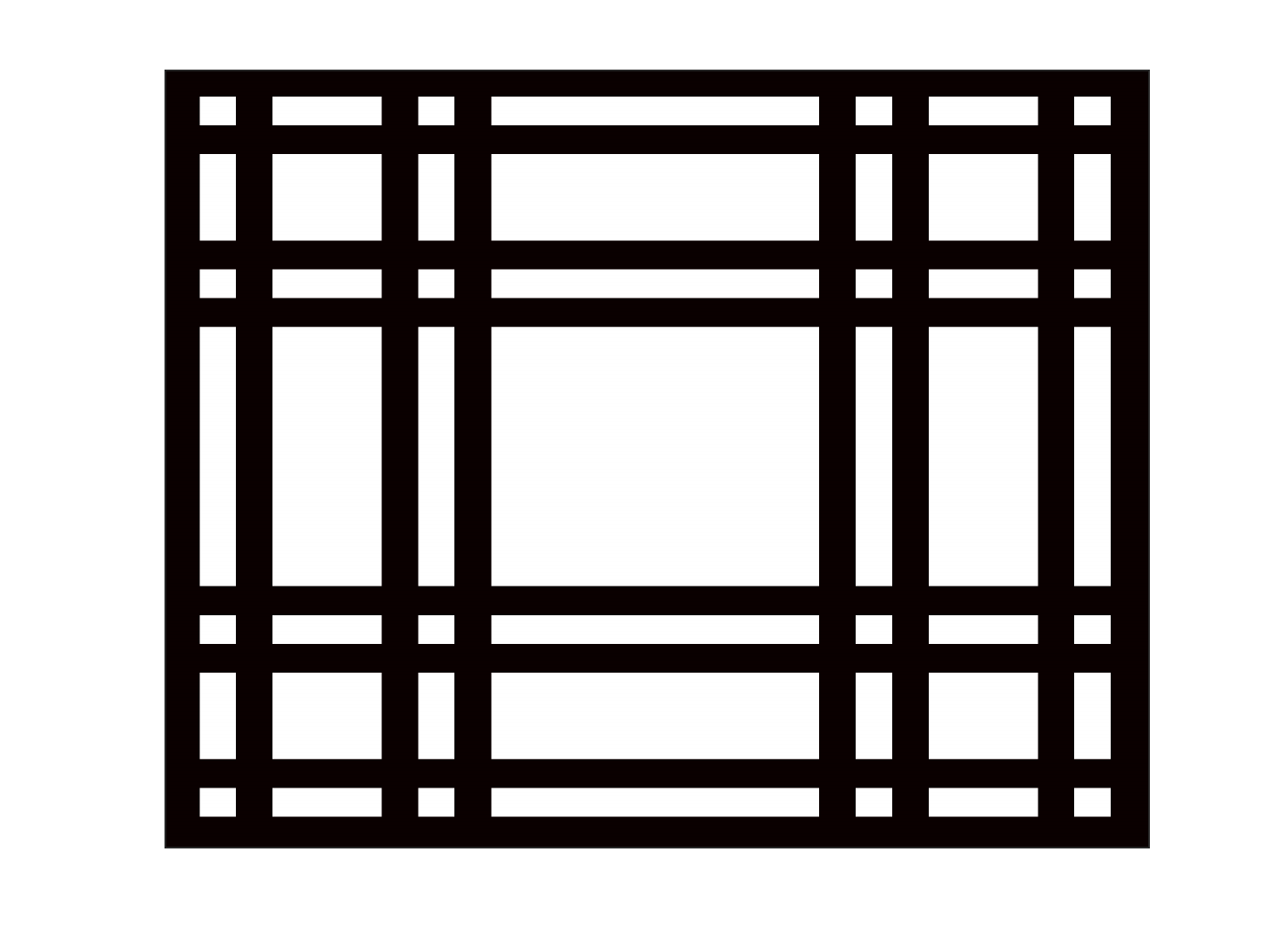}\label{}}
	\subfigure[]{\includegraphics[width = 1.8in]{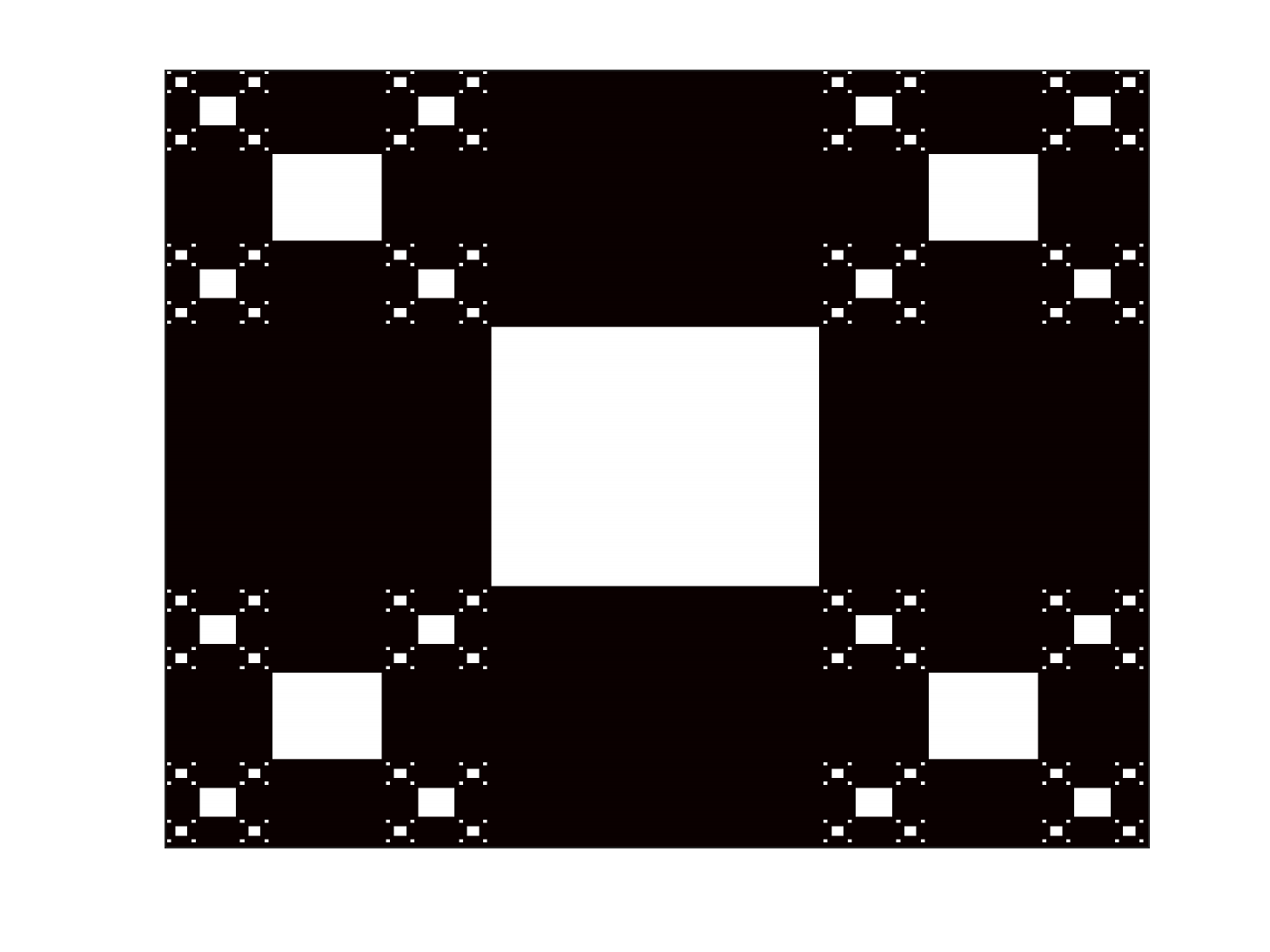}\label{}}
	\caption{Approximations of planar sets generated by (\ref{2DTentMap}) with different conditions of grouping the points}
	\label{CD,pScs}   				
\end{figure}

\begin{figure}[]
	\centering
	\includegraphics[width=0.5\linewidth]{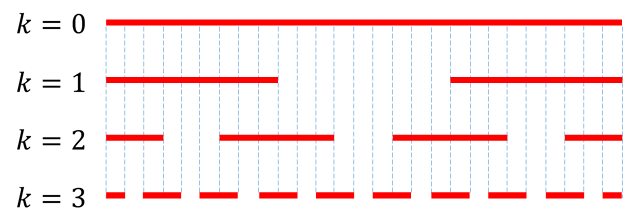}
	\caption{Perforation  set}
	\label{PercolationSet}	
\end{figure}

The construction of the Sierpinski carpet cannot possibly be performed through any arrangement of two dimensional Cantor set and, therefore, a different strategy should be considered. To this end, we shall use maps that construct sets which are similar to Cantor sets in the generation way but different in structure. A suitable set for generating the Sierpinski carpet started with the initial set $ \mathcal{I}=[0, 1] $. The first iteration involves subdividing $ \mathcal{I} $ into three equal intervals and removing the middle open interval $ (\frac{1}{3}, \frac{2}{3}) $. In the second iteration the middle interval is restored and each of the three intervals are again subdivided into three equal subintervals then we remove the middle open intervals $ (\frac{1}{9}, \frac{2}{9}), \; (\frac{4}{9}, \frac{5}{9}) $, and $ (\frac{7}{9}, \frac{8}{9}) $. We continue in the same manner for the succeeding iterations. Figure \ref{PercolationSet} illustrates the first three stages of construction of the set. The purpose of such sets is to cut out successively smaller parts (holes) in the Sierpinski fractals kind. This is why we call these types of sets ``\textit{perforation sets}".

To construct perforation sets, we use the modified tent map
\begin{equation*} \label{ModTentMap}
F(x) = \left\{ \begin{array}{ll}\vspace{2mm}
3 \, [x (\text{mod} \, 1)], & \; \text{if} \;\; x \leq \frac{1}{2} \; \text{or} \; x > 1, \\ 
3 (1-x), & \; \text{if} \;\; \frac{1}{2} < x \leq 1.
\end{array}
\right.
\end{equation*}
The point $ x \in \mathcal{I} $ is excluded from the $ k $th approximation of the set if its $ k $th iteration $ F^k(x) $ does not belong to $ \mathcal{I} $. For the Sierpinski carpet we use a two dimensional version of the modified tent map defined on the domain $ \mathcal{D}=[0, 1]\times[0, 1] $. We consider the iteration

\begin{equation} \label{ModTentMapIter}
\begin{split}	
& x_{n+1} = \left\{ \begin{array}{ll}\vspace{2mm}
3 \, [x_n (\text{mod} \, 1)], & \; \text{if} \;\; x_n \leq \frac{1}{2} \; \text{or} \; x_n > 1, \\ 
3 (1-x_n), & \; \text{if} \;\; \frac{1}{2} < x_n \leq 1,
\end{array} \right. \\
& y_{n+1} = \left\{ \begin{array}{ll}\vspace{2mm}
3 \, [y_n (\text{mod} \, 1)], & \; \text{if} \;\; y_n \leq \frac{1}{2} \; \text{or} \; y_n > 1, \\ 
3 (1-y_n), & \; \text{if} \;\; \frac{1}{2} < y_n \leq 1.
\end{array} \right.
\end{split}		
\end{equation}

To generate the Sierpinski carpet we exclude any point $ (x_0, y_0) \in \mathcal{D} $ if its iteration  $ (x_n, y_n) $ under (\ref{ModTentMapIter}) escapes from $ \mathcal{D} $ such that $ x_n> 1, \; y_n > 1 $ for some natural number $ n $. Figure \ref{SFC} shows the $ 6 $th approximation of the Sierpinski carpet generated by iteration (\ref{ModTentMapIter}). The iteration looks very similar to the FJI but with different criterion for grouping the points. However, it can be classified under FJI type.

Another scheme can be developed by using a map to generate a sequence for each point in a given domain and then applying a suitable criterion to group the points. For that purpose, let us introduce the map
\begin{equation}
\psi_n(x)= B \sin(A_n x),
\label{SCScheme}
\end{equation}
where $ A_n=\pi a^{n-1} $, $ B=\csc\frac{\pi}{b} $, and $ a, b $ are parameters. The recursive formula is defined as follows:
\begin{equation*}
\begin{split}
&\psi_0(x_0):=x_0, \\
&x_n = \psi_n(x_0), \; n=1,2, ... \, .
\end{split}	
\end{equation*}

To construct the perforation set, we start with the interval $ \mathcal{I}=[0, 1] $, and include in the $ k $th approximation of the set each point $ x_0 \in \mathcal{I} $ that satisfies $ |x_k| \leq 1 $. Thus, for Sierpinski carpet, we use a two dimensional version of the map (\ref{SCScheme}) which can be defined in the form
\begin{equation}
\psi_n(x, y)= \big(B \sin(A_n x), \; B \sin(A_n y) \big).
\label{SCI}
\end{equation}
The procedure here is to determine the image sequence $ (x_n, y_n) $ of each point $ (x_0, y_0) \in \mathcal{D} $, i.e.,
\begin{equation}
(x_n, y_n)=\psi_n (x_0, y_0).
\label{SCFJS}
\end{equation}
If we choose $ \mathcal{D}=[0, 1]\times[0, 1] $, the point $ (x_0, y_0) $ is excluded from the set if the condition
\begin{equation}
|x_n| >1, \; |y_n| >1
\label{ExcluCond}
\end{equation}
is satisfied for some $ n \in \mathbb{N} $.

\begin{figure}[]
	\centering	
	\begin{minipage}{0.3\textwidth}
		\subfigure[$ a=b=4 $]{\hspace{-0.3cm}\includegraphics[width = 1.8in]{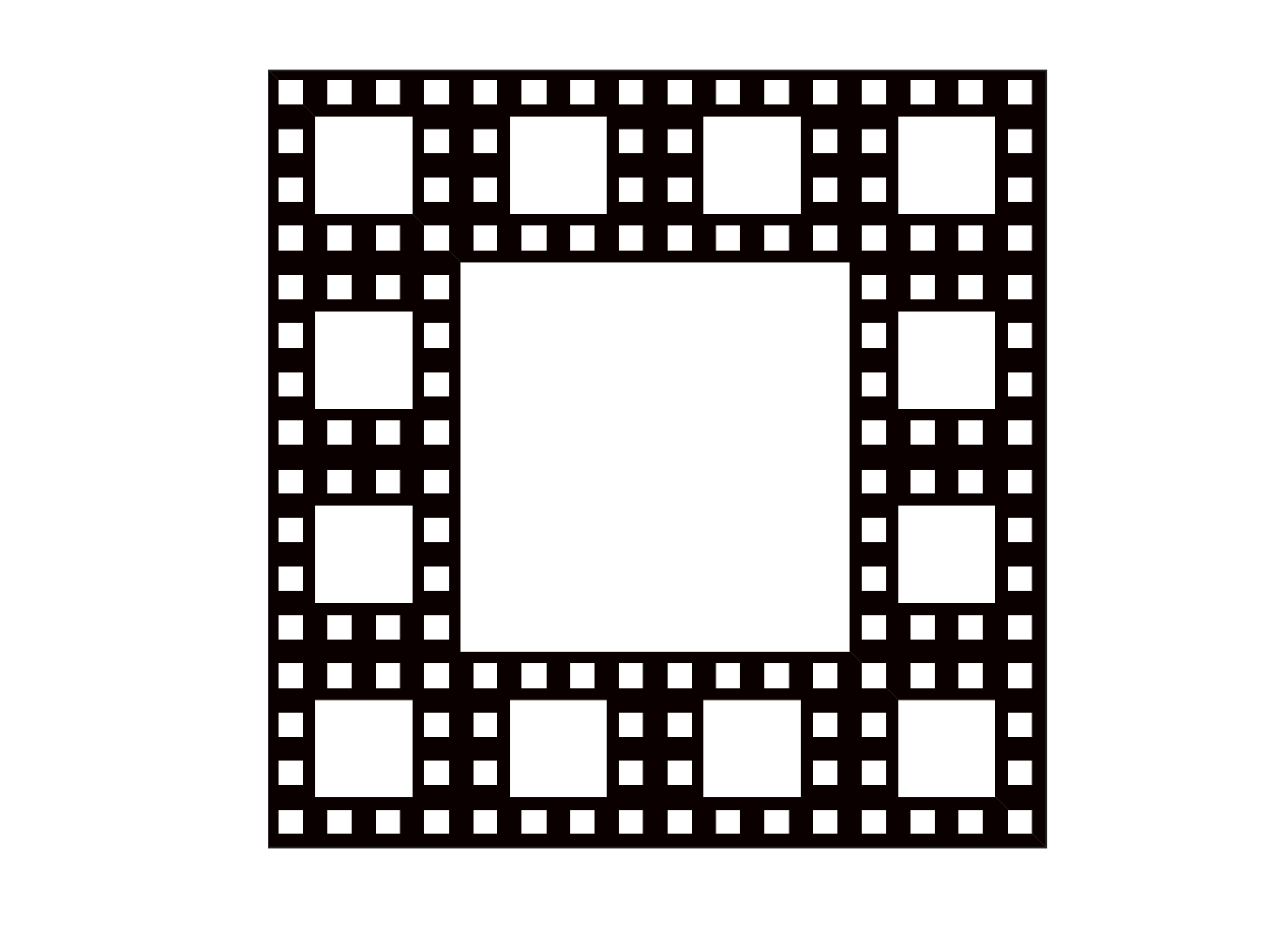}}
		\subfigure[$ a=6, \, b=3 $]{\hspace{-0.3cm}\includegraphics[width = 1.8in]{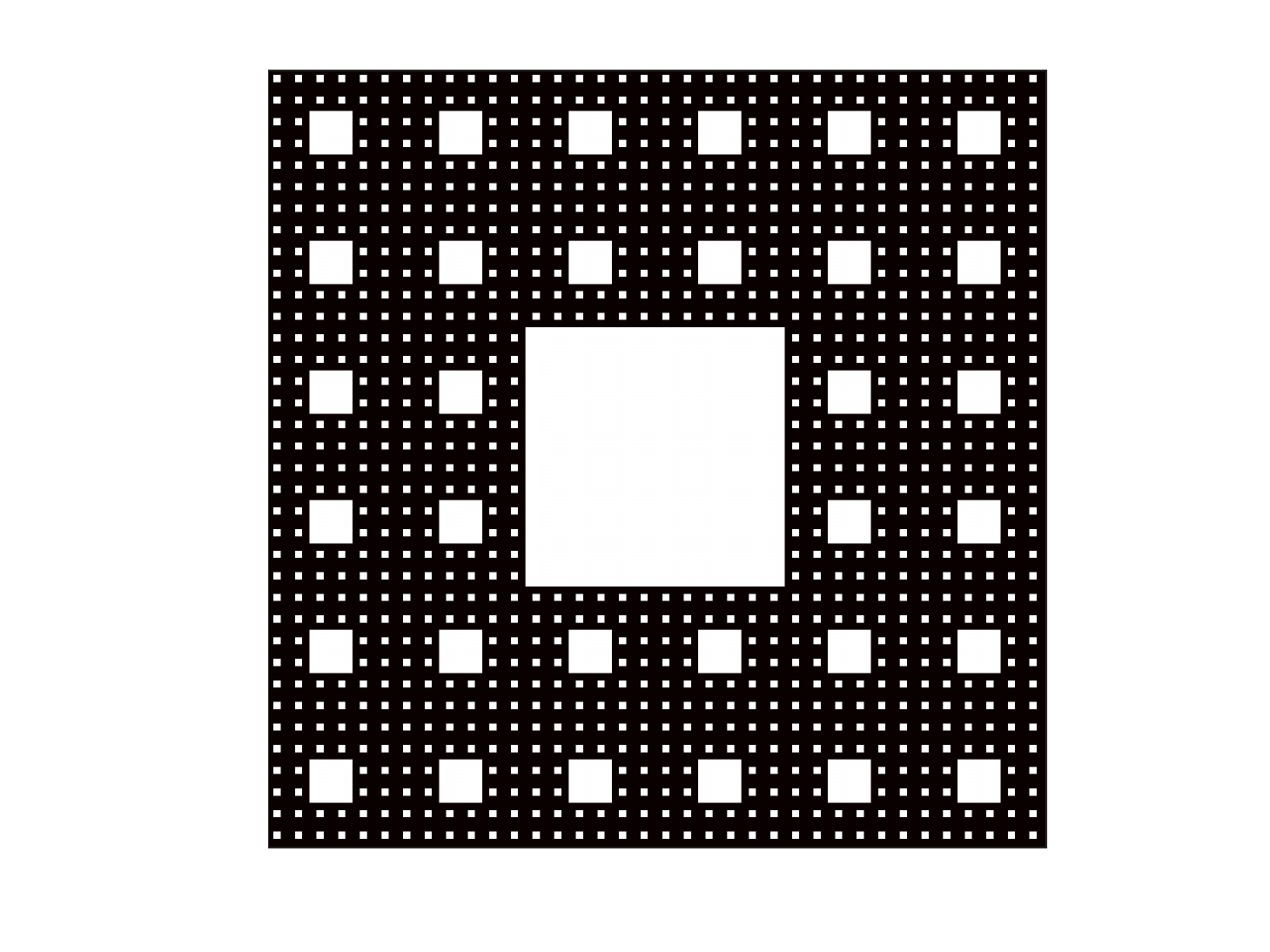}}
	\end{minipage} 
	\begin{minipage}{0.3\textwidth}		
		\subfigure[$ a=3, \, b=4 $]{\hspace{-0.3cm}\includegraphics[width = 1.8in]{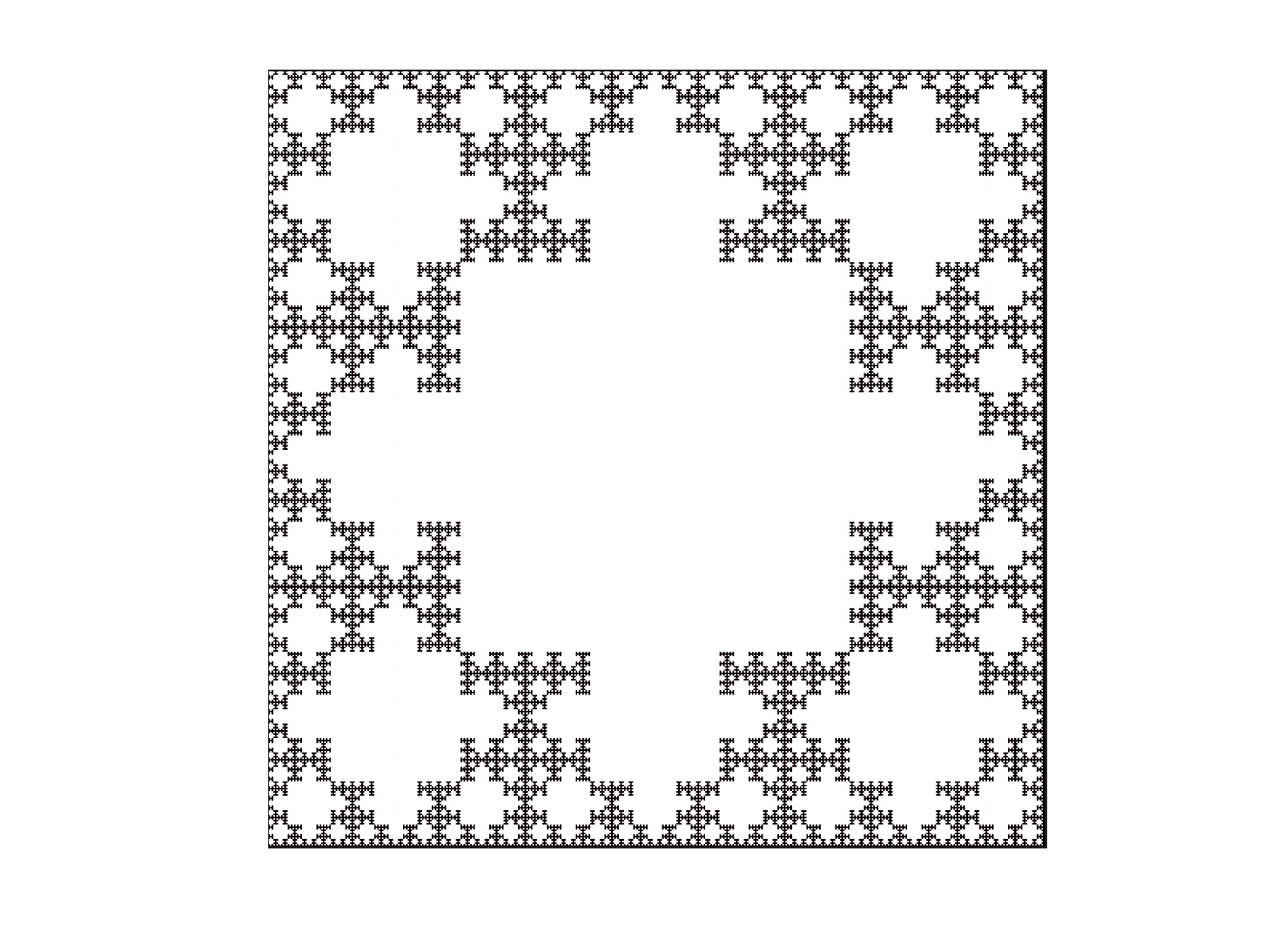}}
		\subfigure[$ a=3, \, b=4 $]{\hspace{-0.3cm}\includegraphics[width = 1.8in]{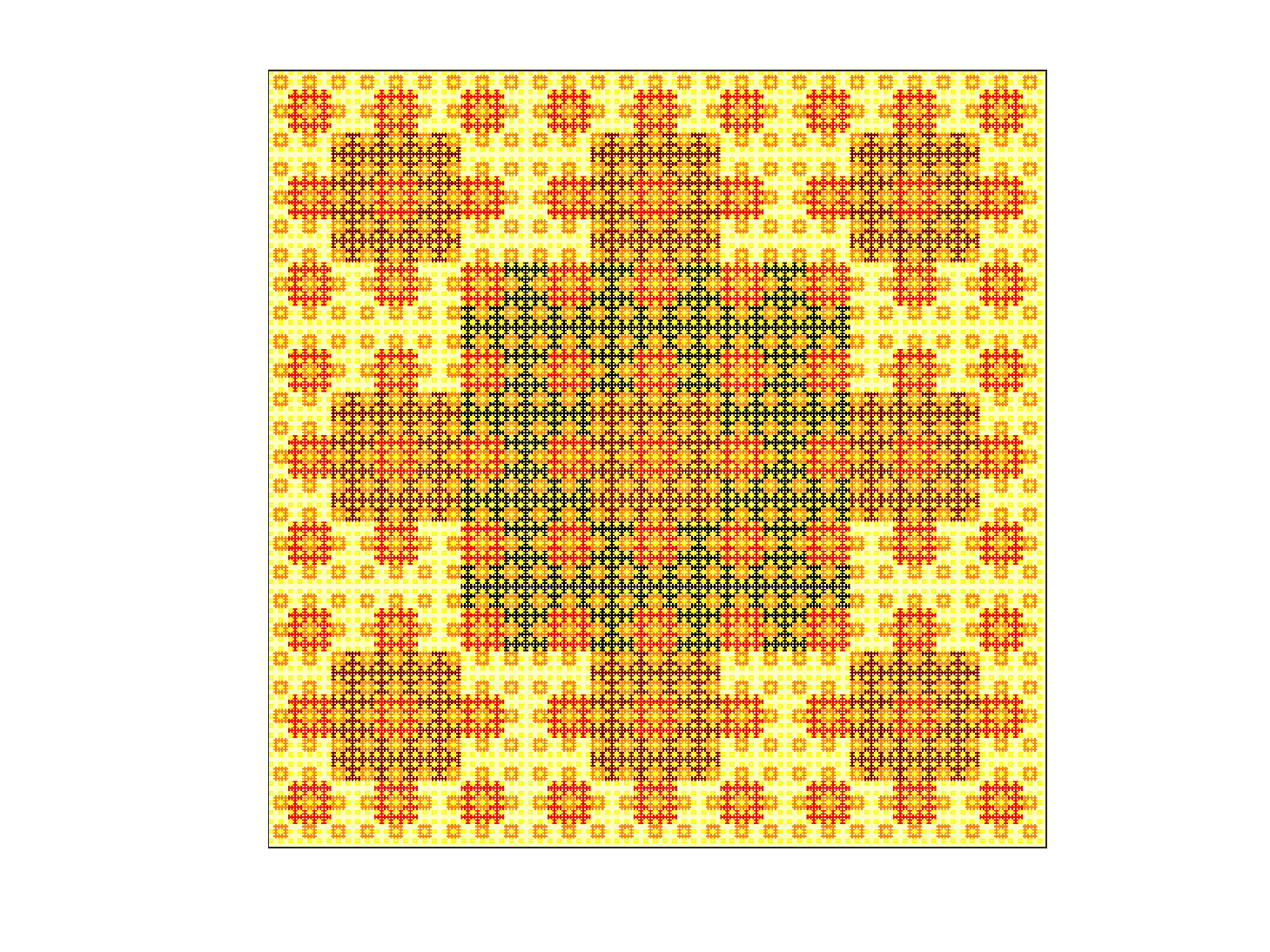}}
	\end{minipage} 
	\begin{minipage}{0.3\textwidth}		
		\subfigure[$ a=2, \, b=3.5 $]{\hspace{-0.3cm}\includegraphics[width = 1.8in]{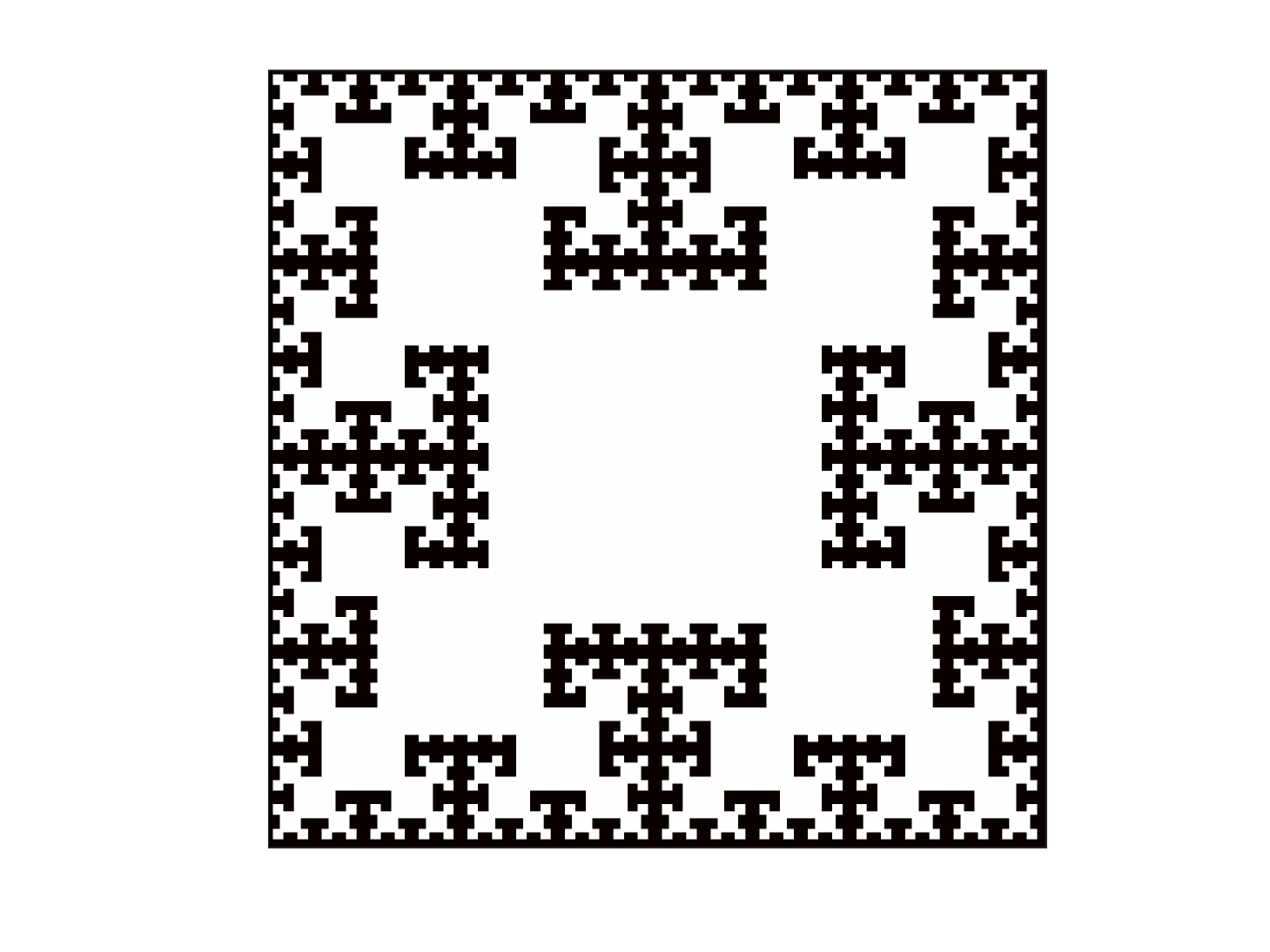}}
		\subfigure[$ a=2, \, b=3.5 $]{\hspace{-0.3cm}\includegraphics[width = 1.8in]{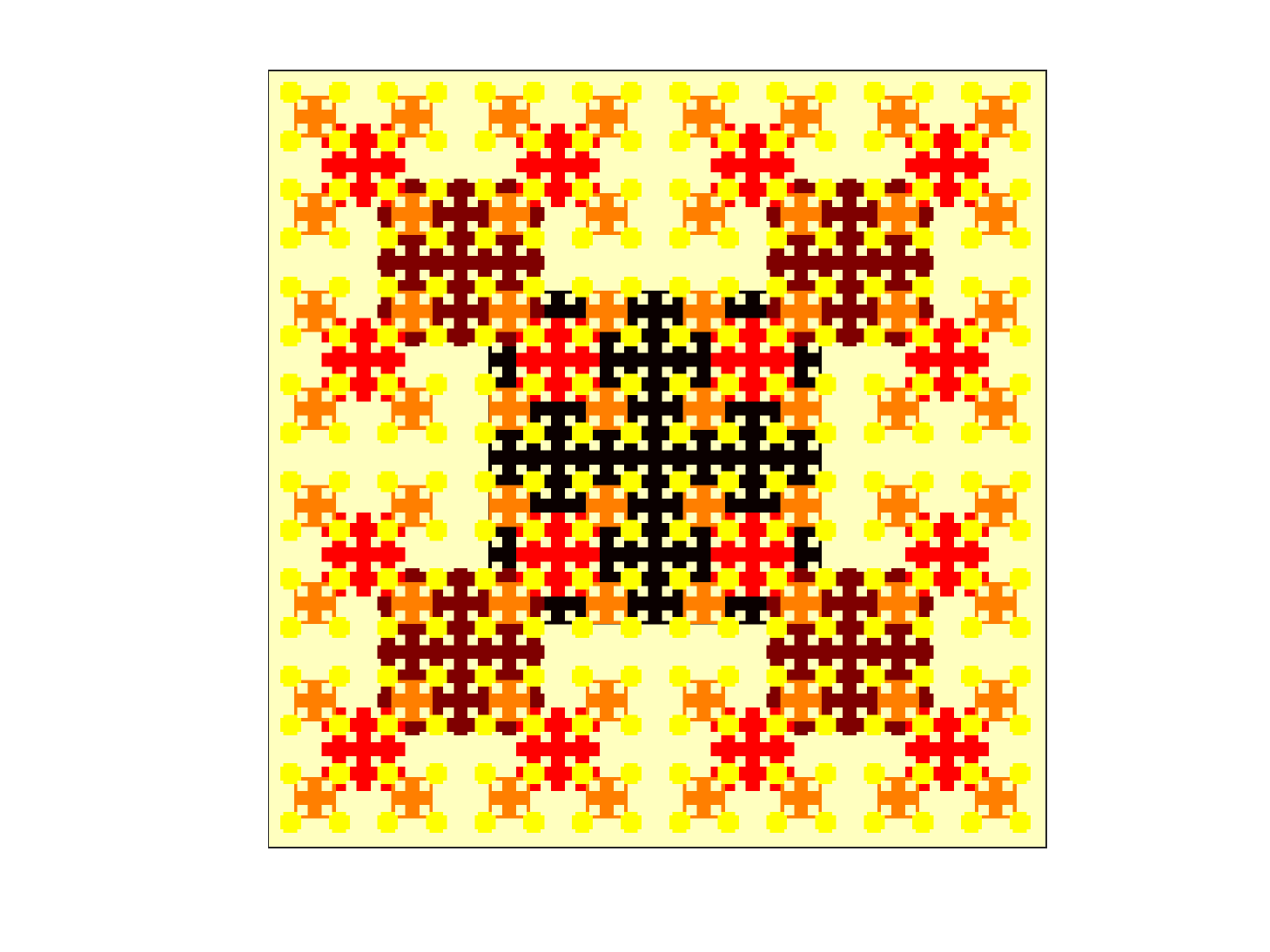}}
	\end{minipage} %}}
	\caption{Sierpinski Carpets}
	\label{SLikeCarpet}   				
\end{figure}

For the values of the parameters $ a=3 $ and $ b=3 $, the scheme gives the classical Sierpinski carpet (the simulation result for the $ 6 $th approximation is identical to Fig. \ref{SFC}). Figure \ref{SLikeCarpet} shows other carpets generated by (\ref{SCI}) with different values of the parameters $ a $ and $ b $, and for limited stages. The colors that appear in the parts (d) and (f) of the figure are related to the sequences generated by (\ref{SCFJS}) such that the color of each point in the carpets depends on the smallest number of stages $ n $ that satisfies condition (\ref{ExcluCond}).

This scheme is quite different than the usual procedure of FJI since iterations are not utilized and a different criterion is applied for grouping the points. However, we shall see later that the idea of the FMI can be applied for this type which allows to map and introduce dynamics for the constructed carpets.

A more similar iteration to that of Fatou and Julia can be constructed by finding the inverse $ x=\psi_n^{-1}(x_n) $ in (\ref{SCScheme}) and then substituting in $ x_{n+1}=\psi_{n+1}(x) $ to get the autonomous iteration
\begin{equation}
x_{n+1}= B \sin \Big( a \sin^{-1}\frac{x_n}{B} \Big).
\label{AutoSCI}
\end{equation}

\begin{figure}[]
	\centering
	\subfigure[$ a=3, \; b=3 $]{\includegraphics[width = 2.0in]{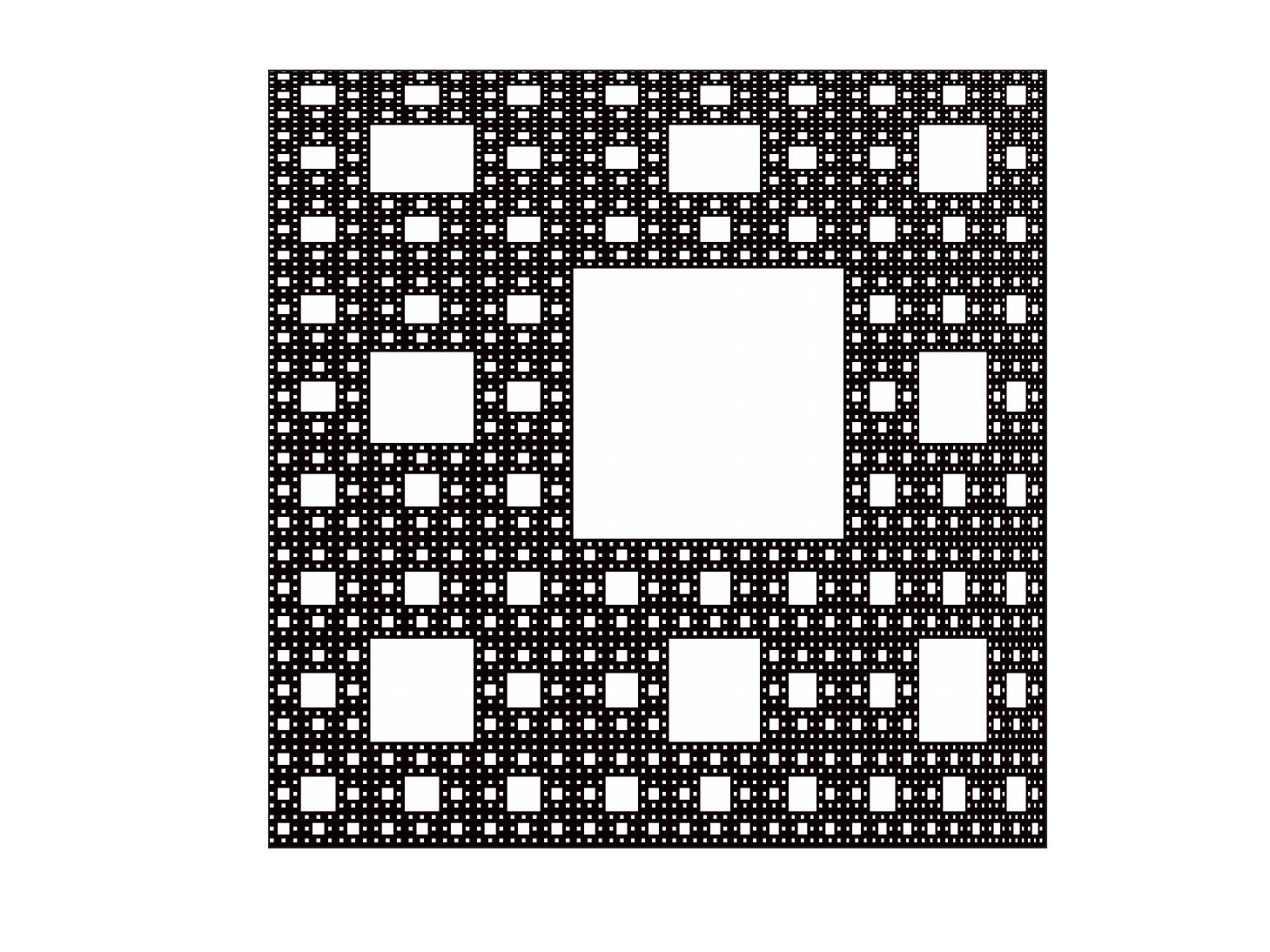}\label{}} 
	\subfigure[$ a=4, \; b=4 $]{\includegraphics[width = 2.0in]{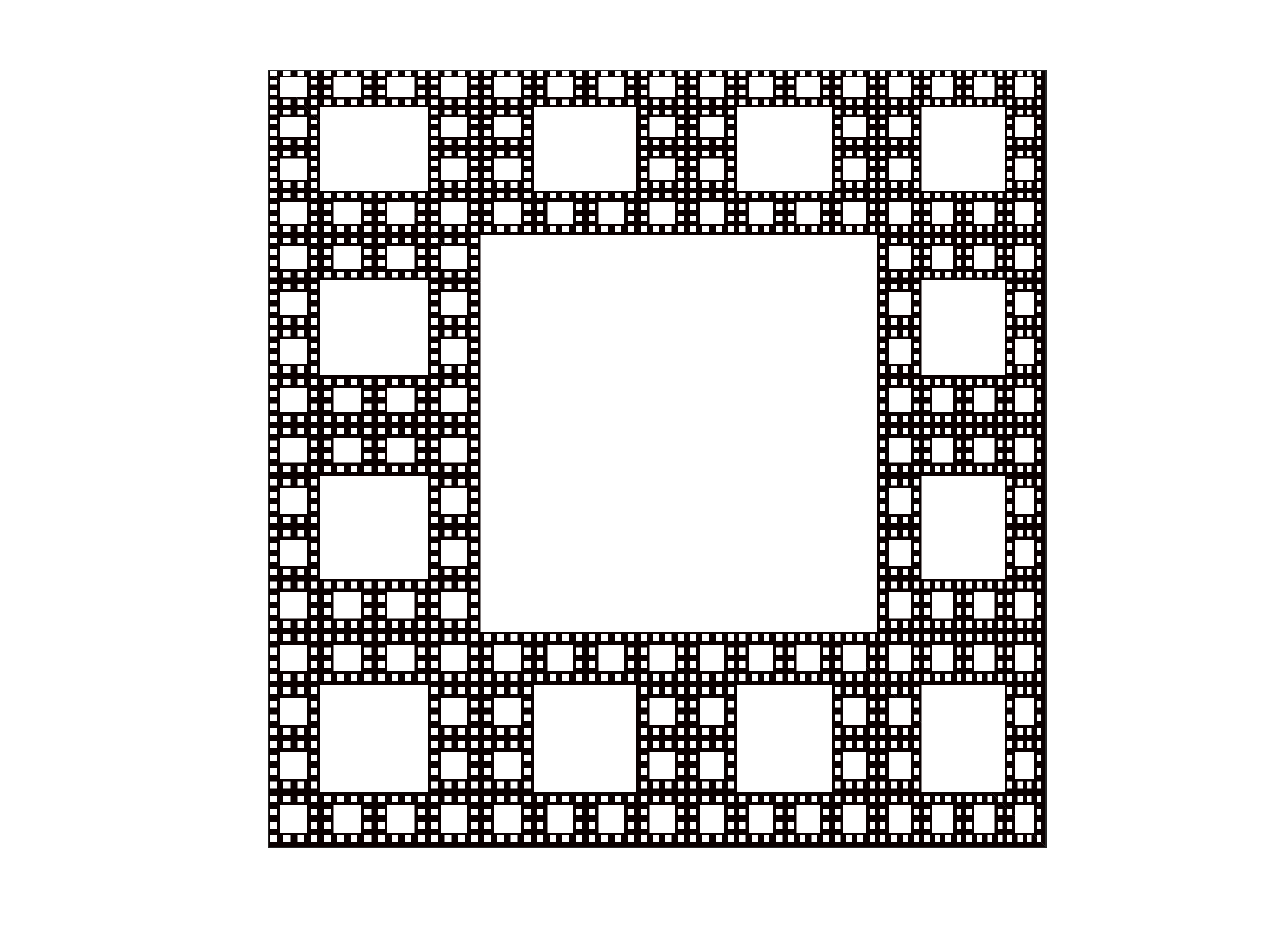}\label{}}
	\caption{Sierpinski carpets by the FJI (\ref{AutoSCI})}
	\label{IrregSC}   				
\end{figure}

The carpets shown in Fig. \ref{IrregSC} are simulated by using a two dimensional form of iteration (\ref{AutoSCI}), and they are irregular types of  fractals. These fractals have asymmetric similarities and they are not categorized under random fractals \cite{Feldman}.

\subsection{\normalsize Sierpinski gasket}

\begin{figure}[]
	\centering
	\begin{minipage}{0.45\textwidth}
		\centering
		\subfigure[Coordiate system for Sierpinski gasket construction]{\includegraphics[width = 2.0in]{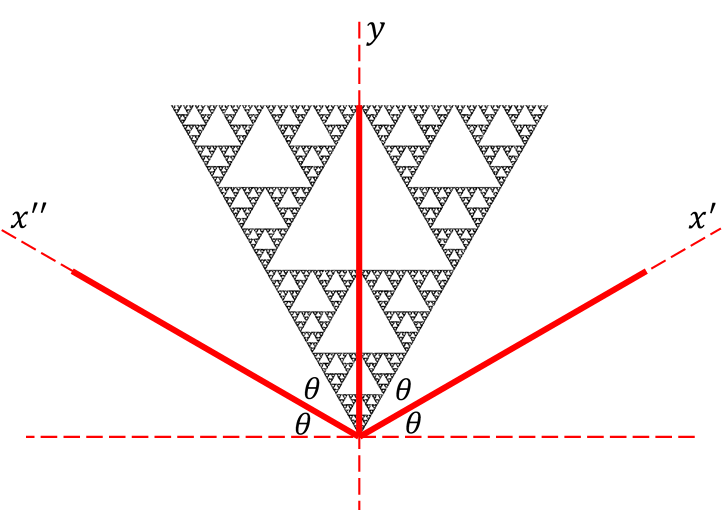}\label{SGCoord}}
	\end{minipage}
	\begin{minipage}{0.45\textwidth}
		\centering
		\subfigure[Perforation set]{\includegraphics[width = 1.8in]{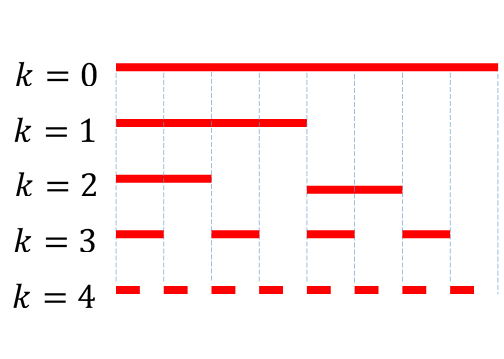}\label{PercolationSet2}}
	\end{minipage}
	\caption{Sierpinski gasket construction}
	\label{SGConstruc}   				
\end{figure}

For constructing the Sierpinski gasket we again use the perforation sets, and in this case we introduce a special coordinate system shown in Fig. \ref{SGCoord}. The system consists of three non-rectangular plane axes denoted by $ x' $, $ x'' $, and $ y $, and the thick red lines in the figure represent the perforation set constructed in Fig. \ref{PercolationSet2}.

We start with the initial set $ \mathcal{D} $ which is a unit equilateral triangle defined by $ \mathcal{D}=\{ (x, y) \in \mathbb{R}^2 : \frac{-1}{\sqrt{3}}y \leq x \leq \frac{1}{\sqrt{3}}y , \; 0 \leq y \leq \frac{\sqrt{3}}{2} \} $. For each point $ (x, y) \in \mathcal{D} $, we detect the triple $ (x', x'', y) $, where $ x' $ and $ x'' $ are the projections of $ (x, y) $ on the $ x' $ and $ x'' $ axes, respectively. To examine whether the point $ (x, y) $ belongs to the gasket, we set up the recursive formula
\begin{equation}
(x'_n, x''_n, y_n)= \Big( \alpha \big(A_n x' \big) , \; \beta \big(A_n x'' \big) , \; \gamma \big(A_n y \big) \Big),
\label{SGI}
\end{equation}
where $ \alpha, \beta $ and $ \gamma $ are functions, $ A_n=\frac{2}{\sqrt{3}} \pi a^n $, and $ a $ is a parameter. The point $ (x, y) $ is excluded from the resulting gasket if $ x'_n > 0, \, x''_n > 0 $ and $  y_n < 0 $ for some $ n \in \mathbb{N} $.

\begin{figure}[]
	\centering	
	\begin{minipage}{0.9\textwidth}
		\centering	
		\subfigure[$ \alpha(x)= \beta(x)= \gamma(x)=\sin(x), \; a=4 $]{\includegraphics[width = 1.8in]{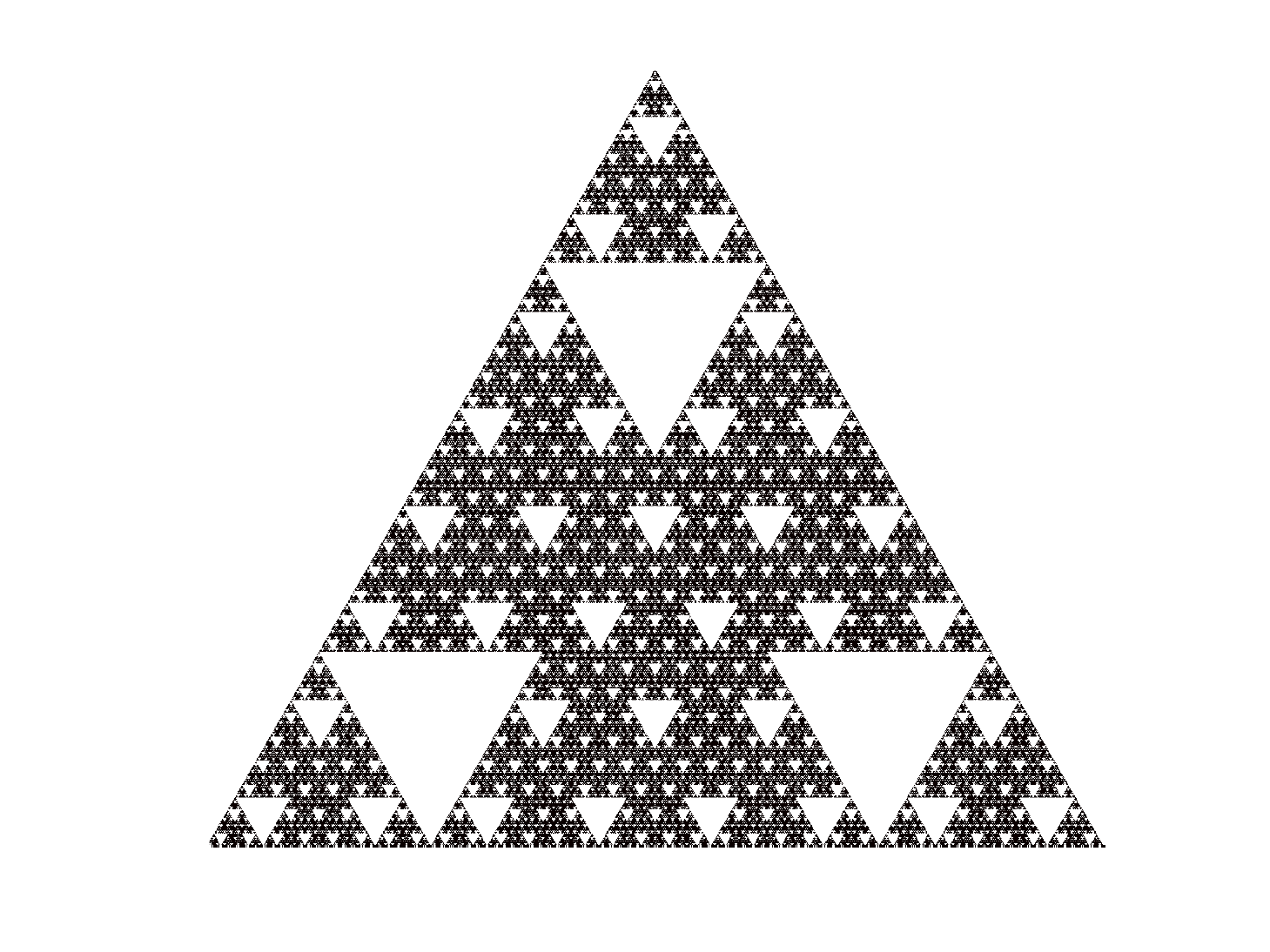}} \hspace{1.5cm}
		\subfigure[$ \alpha(x)= \beta(x)= \gamma(x)=\cos(x), \; a=2 $]{\includegraphics[width = 1.8in]{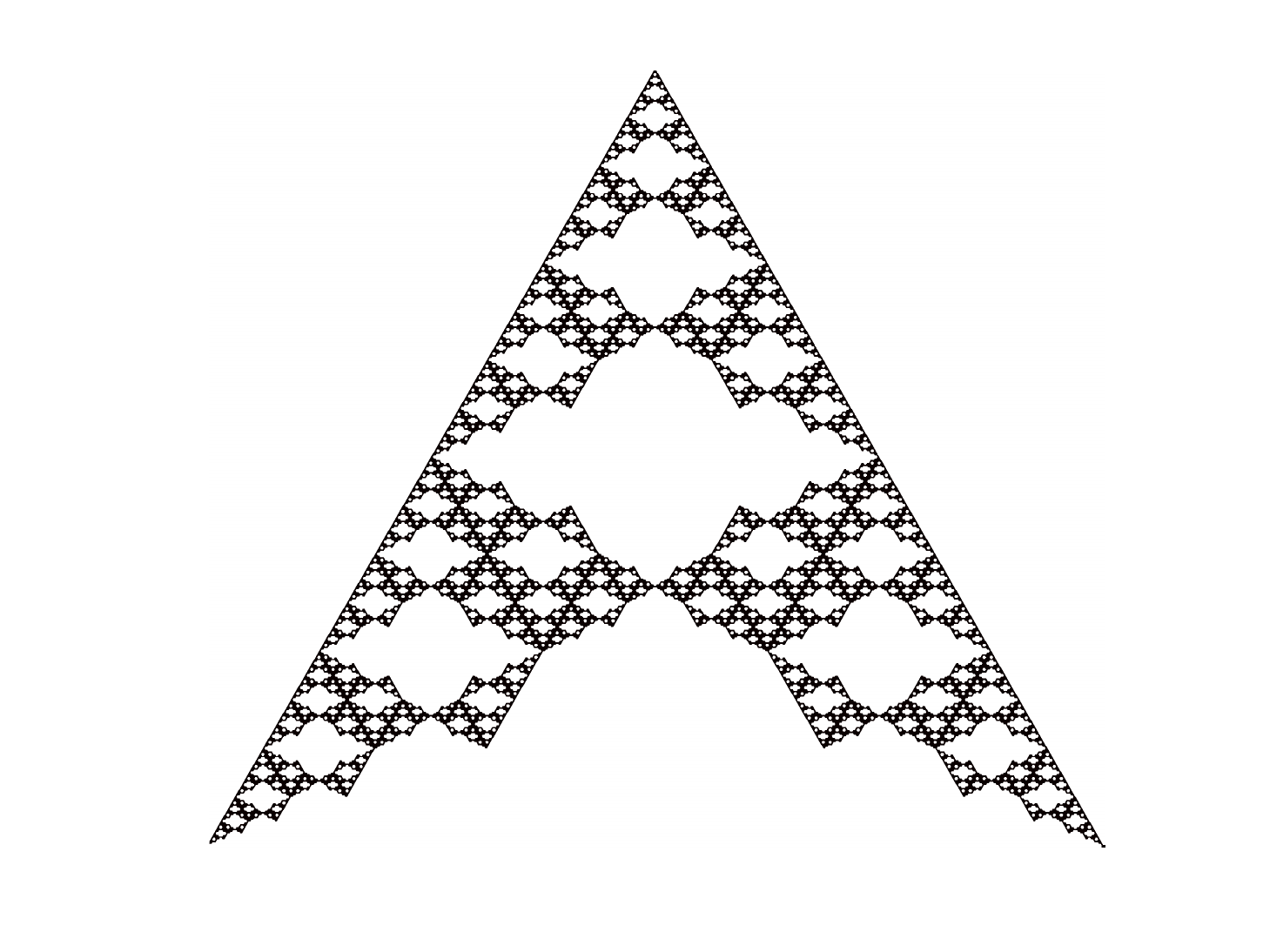}}
	\end{minipage} 
	\begin{minipage}{0.9\textwidth}	
		\centering		
		\subfigure[$ \alpha(x)=\tan(x), \; \beta(x)= \gamma(x)=\cos(x), \; a=2 $]{\includegraphics[width = 1.8in]{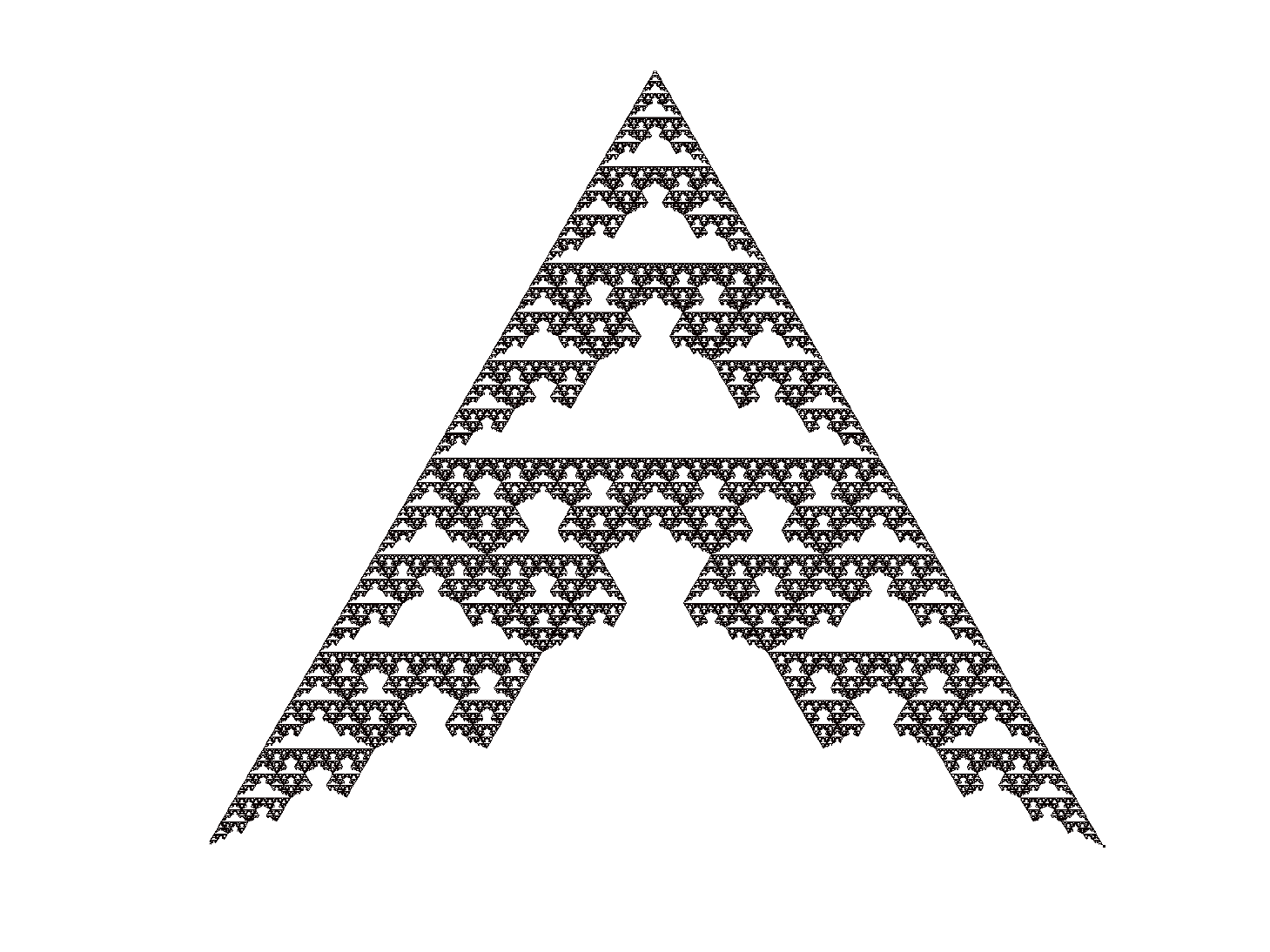}}\hspace{1.5cm}
		\subfigure[$ \alpha(x)=\tan^{-1}(x), \; \beta(x)= \gamma(x)=\cos(x), \; a=7 $]{\includegraphics[width = 1.8in]{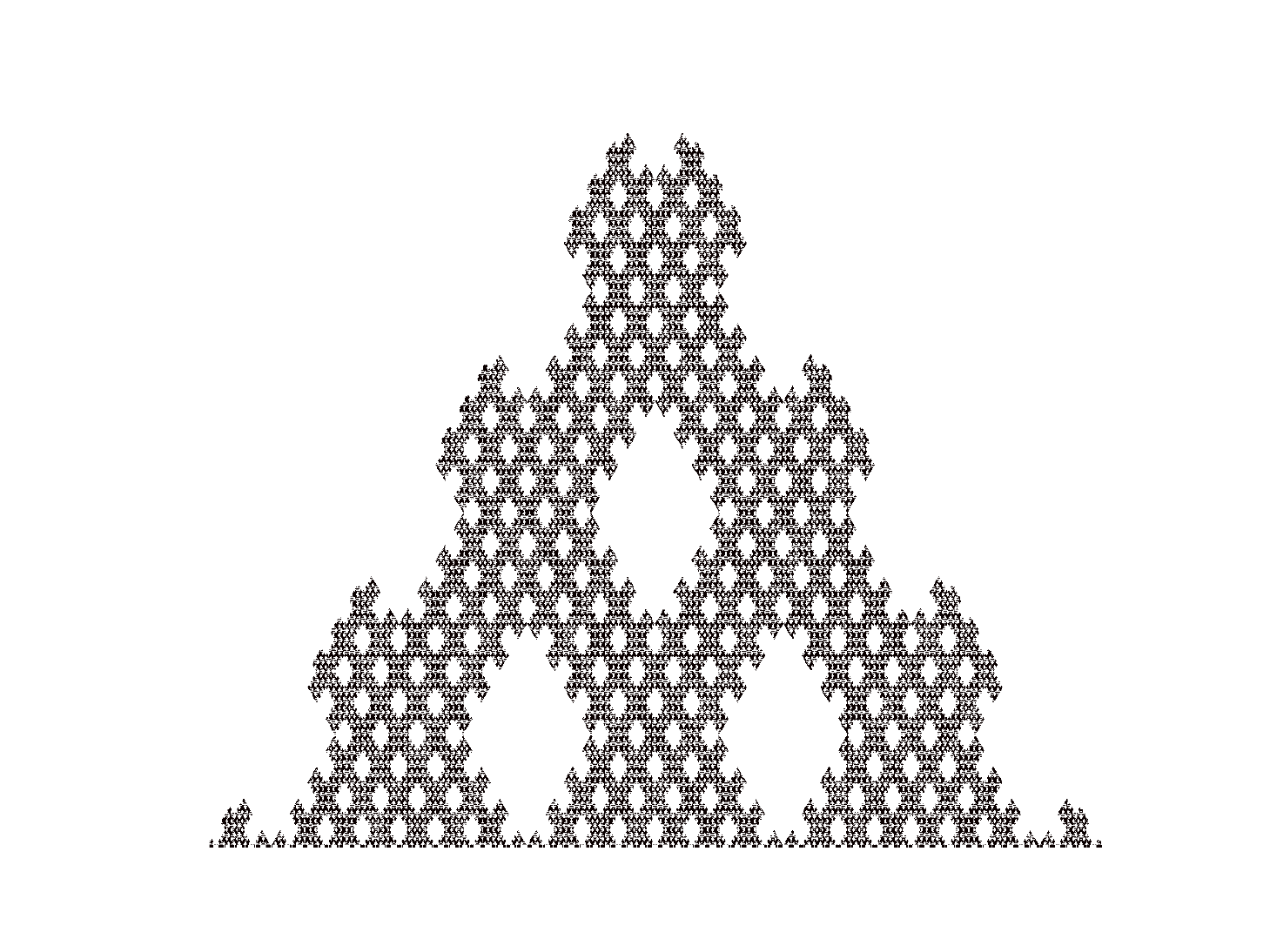}}
	\end{minipage} 
	\caption{Sierpinski Gaskets}
	\label{SLikeGasket}   				
\end{figure}

For $ \alpha(x)= \beta(x)= \gamma(x)= \sin x $ and $ a =2 $, the resulting gasket is the classical Sierpinski gasket. Figure \ref{SFG} shows the $ 8 $th approximation of the Sierpinski gasket generated by iteration (\ref{SGI}). Examples of other gaskets with different choices of the functions $ \alpha, \beta, \gamma $ and the parameter $ a $ are shown in Fig. \ref{SLikeGasket}.

\section{\normalfont MAPPINGS}		

In this part of the paper we give procedures for mapping the Sierpinski fractals through the schemes introduced in the preceding section. 

\subsection{\normalsize Mapping of carpets}

To map the carpets generated by the scheme (\ref{SCI}), we use the idea of FMI. Let $ \varPhi :  \mathcal{D} \to  \mathcal{D'} $ be an invertible function defined by
\begin{equation}
\varPhi(x, y)= (\phi_1,\phi_2)(x, y),
\label{PhiMap}
\end{equation}
with the inverse
\begin{equation}
\varPhi^{-1}(\xi, \eta )= (\phi_3,\phi_4)(\xi, \eta ).
\label{InvPhiMap}
\end{equation}
Then the fractal mapping scheme can be defined as
\begin{equation} \label{SierpinsliFMI}
\varPhi^{-1}(\xi_n, \eta_n)= \psi_n \big( \varPhi^{-1}(\xi_0, \eta_0) \big).
\end{equation}
This scheme transforms a carpet $ \mathcal{F} $ into a new carpet $ \mathcal{F}_\varPhi $, and the following theorem shows that the set $ \mathcal{F}_\varPhi $ is merely the image of $ \mathcal{F} $ under the map $ \varPhi $.

\begin{theorem} \label{Thm3}
	$ \mathcal{F}_\varPhi=\varPhi(\mathcal{F}) $.
\end{theorem}
\noindent {\bf Proof.}
First, we show that $ \mathcal{F}_\varPhi \subseteq \varPhi(\mathcal{F}) $. Let $ (\xi, \eta) \in \mathcal{F}_\varPhi $, which means that for formula (\ref{SierpinsliFMI}), if $ (\xi_0, \eta_0) = (\xi, \eta) $, then at least one of $ |u_n| $ and $ |v_n| $ is less than or equal to $ 1 $ for all $  n \in \mathbb{N}$, where $ (u_n, v_n)= \varPhi^{-1}(\xi_n, \eta_n) $. This implies that $ (u, v)= (u_0, v_0) \in \mathcal{F} $. Thus $ (\xi, \eta) \in \varPhi(\mathcal{F}) $.

For the reverse inclusion, suppose that $ (\xi, \eta) \in \varPhi(\mathcal{F}) $, i.e., there exists $ (x, y) \in \mathcal{F} $ such that $ \varPhi(x, y)= (\xi, \eta) $ and $ (x_0, y_0)=(x, y) $ with formula (\ref{SCFJS}) in which at least one of $ |x_n| $ and $ |y_n| $ is less than or equal to $ 1 $ for all $  n \in \mathbb{N}$. This directly implies that the sequence $ \varPhi^{-1}(\xi_n, \eta_n) = (x_n, y_n) $ satisfies (\ref{SierpinsliFMI}) and $ (\xi, \eta) \in \mathcal{F}_\varPhi $. $ \square $

The following question arises here: Is the mapped carpet a fractal? The answer is ``\textit{yes}" if the map $ \varPhi $ satisfies a bi-Lipschitz condition. This result is stated in the following lemma.

\begin{lemma} \label{Lem1} \cite{Falconer}
	Let $ E \subseteq \mathbb{R}^n $. If $ f: E \rightarrow \mathbb{R}^m $ is a bi-Lipschitz function, i.e., there exist real numbers $ l_1, l_2 > 0 $ such that
	\begin{equation*}  \label{Bi-Lip cond}
	l_1 |u-v| \leq |f(u)-f(v)| \leq l_2 |u-v|,
	\end{equation*}
	for all $ u, v \in E $, then 
	\[ \dim_H \big(f(E)\big)  = \dim_H(E), \]
	where $ \dim_H $ denotes the Hausdorff  dimension.
\end{lemma}

\begin{figure}[]
	\centering
	\subfigure[]{\includegraphics[width = 2.2in]{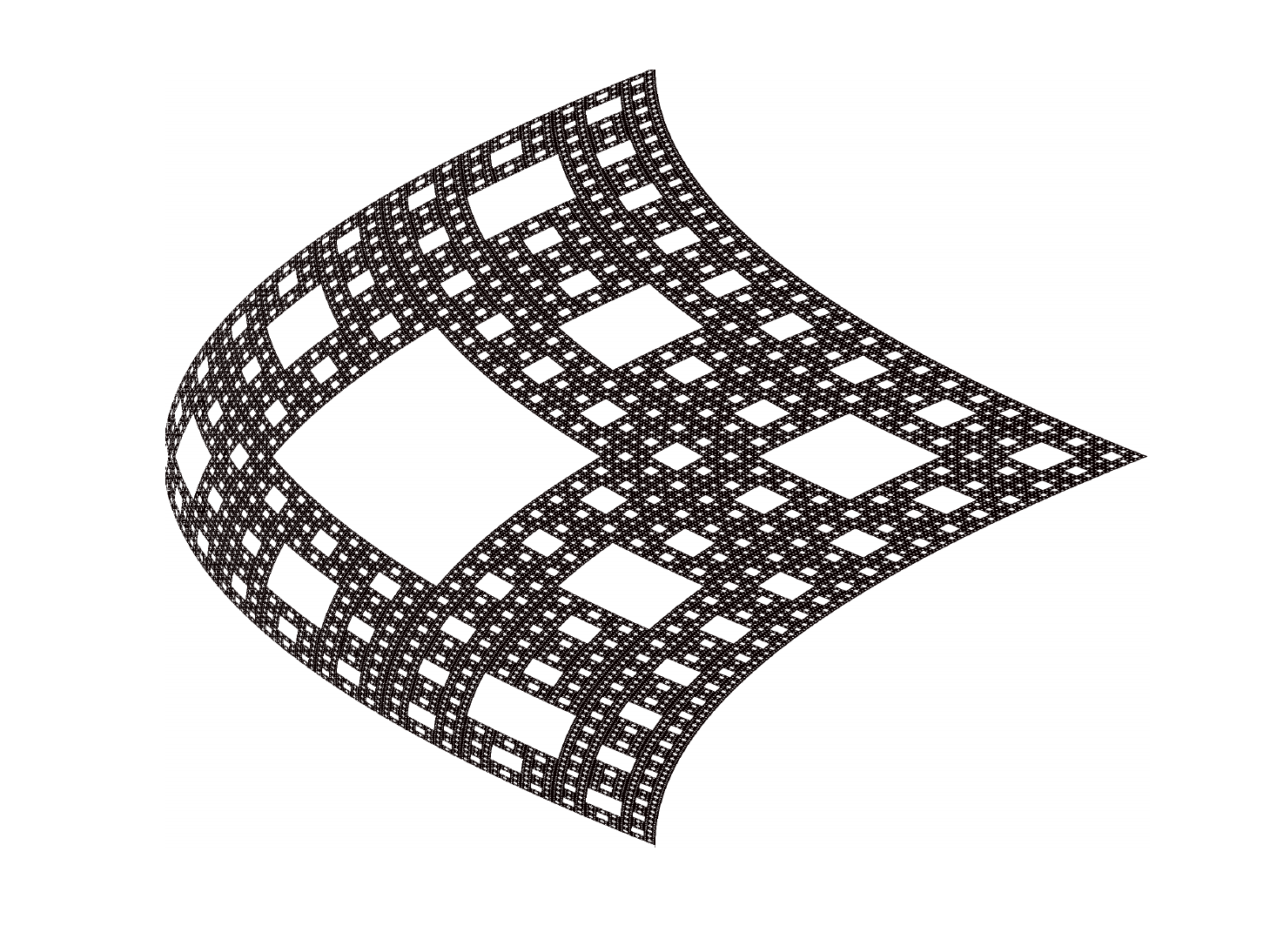}\label{}} 
	\subfigure[]{\includegraphics[width = 2.2in]{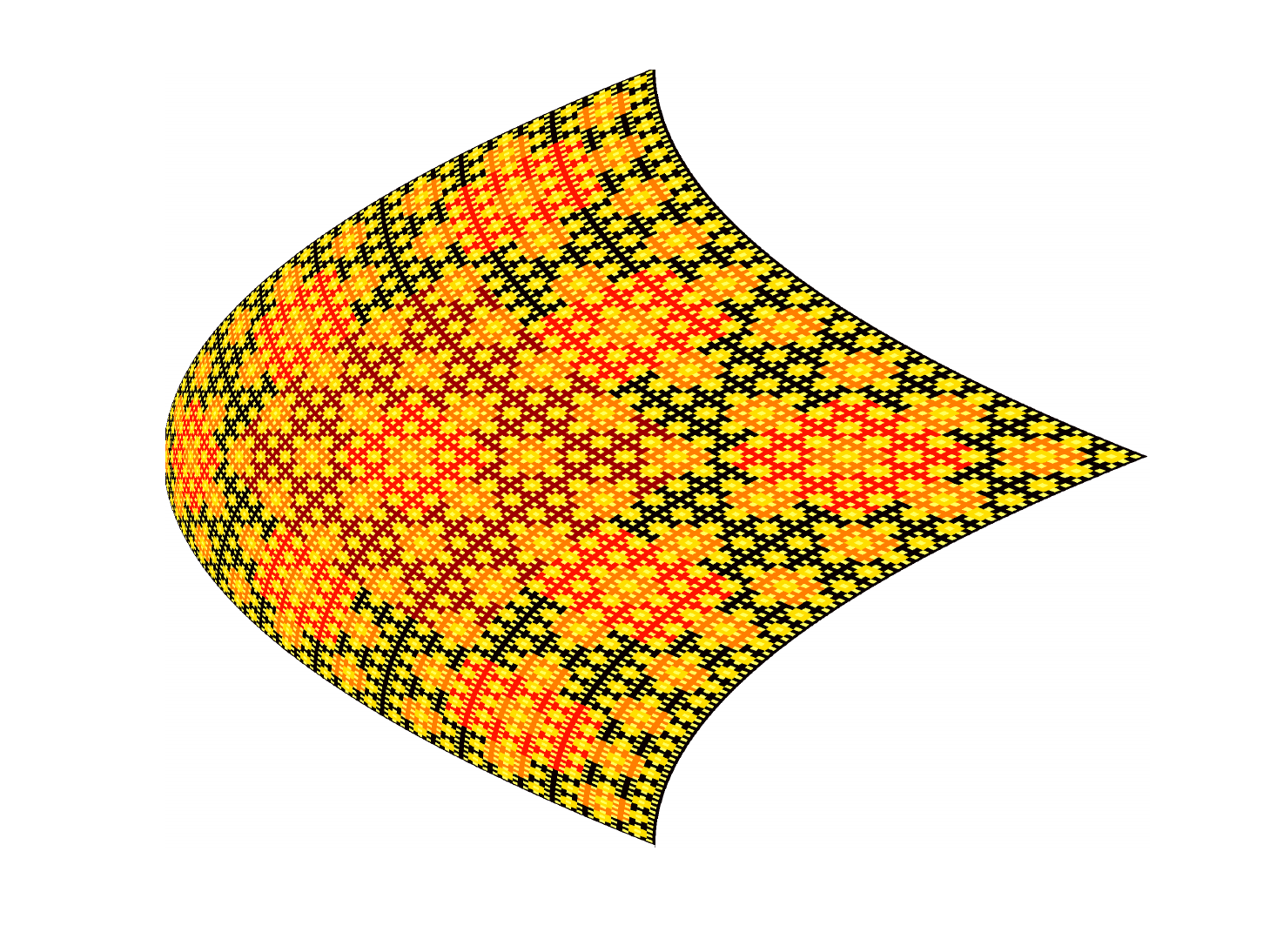}\label{}}
	\caption{The images of carpets with (a): $ a=b=3 $, (b): $ a=3, \, b=4 $}
	\label{MapSC1}   				
\end{figure}

\begin{figure}[]
	\centering
	\subfigure[]{\includegraphics[width = 2.2in]{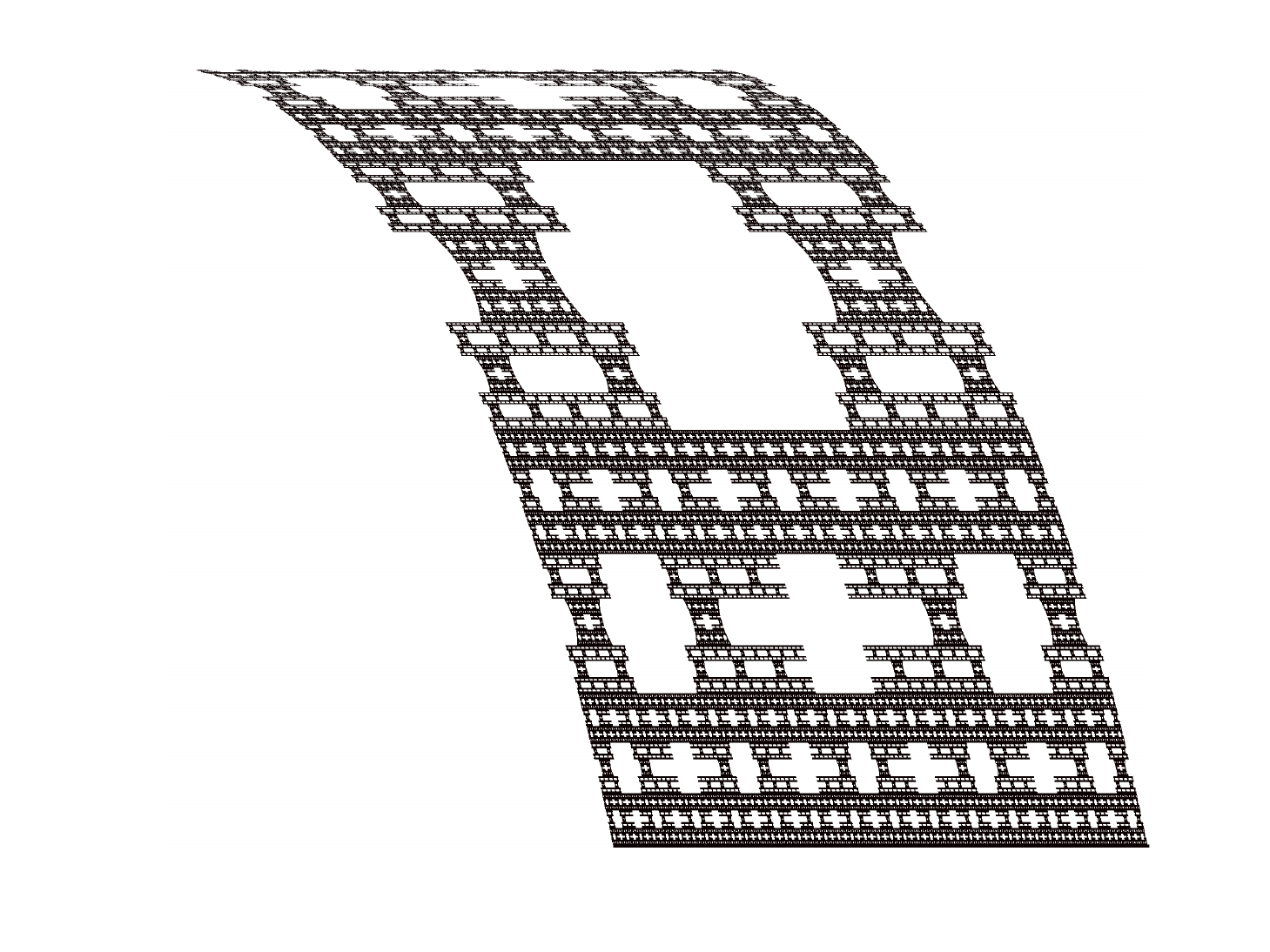}\label{}} 
	\subfigure[]{\includegraphics[width = 2.2in]{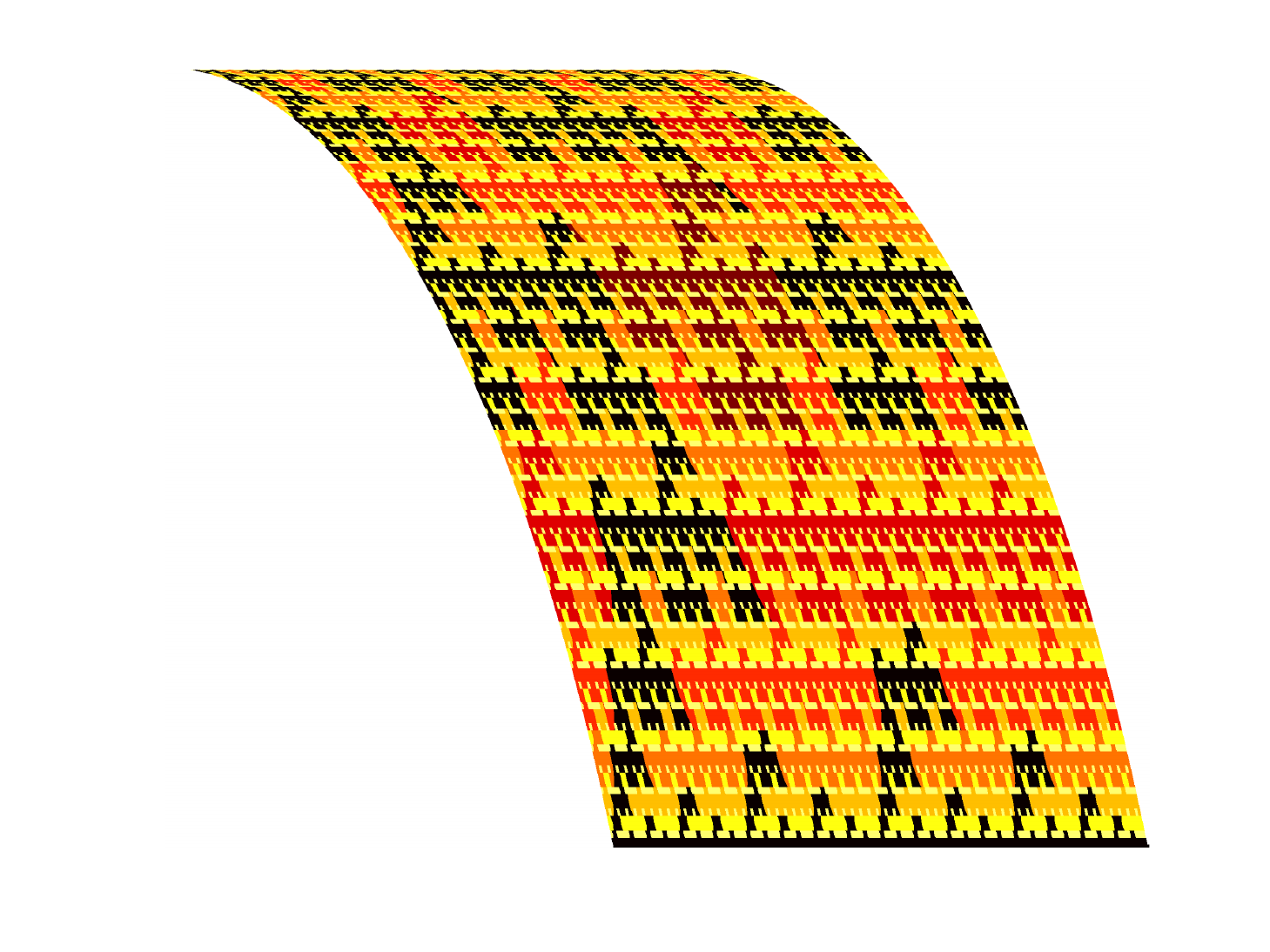}\label{}}
	\caption{The image of carpet with (a): $ a=2, \, b=3.5 $, (b): $ a=2, \, b=1.5 $}
	\label{MapSC2}   				
\end{figure}

For our next examples, we shall use bi-Lipschitz functions to ensure that the mapped carpets are fractals. In order to obtain the mapped Sierpinski carpet $ \mathcal{F}_\varPhi $, we restrict the domain of (\ref{SierpinsliFMI}) only to the points $ (\xi, \eta) $ that belong to the mapped domain $ \mathcal{D'} $, thus Eq. (\ref{SierpinsliFMI}) becomes
\[ (\xi_n, \eta_n)= \varPhi \big(\psi_n (x_0, y_0) \big), \quad (x_0, y_0) \in \mathcal{D}. \]
More precisely, by using Eqs. (\ref{SCI}) and (\ref{InvPhiMap}) in (\ref{SierpinsliFMI}), we have
\[ \varPhi^{-1}(\xi_n, \eta_n)= \Bigg(\frac{\sin(a^{n-1} \pi \phi_3(\xi_0, \eta_0))}{\sin(\frac{\pi}{b})} , \; \frac{\sin(a^{n-1} \pi \phi_4(\xi_0, \eta_0))}{\sin(\frac{\pi}{b})} \Bigg).	\]
Now letting $ \frac{\sin(a^{n-1} \pi \phi_3(\xi_0, \eta_0))}{\sin(\frac{\pi}{b})} = X_n $ and $ \frac{\sin(a^{n-1} \pi \phi_4(\xi_0, \eta_0))}{\sin(\frac{\pi}{b})} = Y_n $ and using Eq. (\ref{PhiMap}), we get	
\begin{equation} \label{MapCarpetsEq}
(\xi_n, \eta_n)= \Big( \phi_1(X_n, Y_n), \; \phi_2(X_n, Y_n) \Big).
\end{equation}
The semi-iteration (\ref{MapCarpetsEq}) is applied for each point $ (\xi_0, \eta_0) \in \mathcal{D'} $, and the point is excluded from the image set $ \mathcal{F}_\varPhi $ if $ |\phi_3(\xi_n,\eta_n)| >1, \; |\phi_4(\xi_n,\eta_n)| >1 $ for some $ n \in \mathbb{N} $. Figures \ref{MapSC1} and \ref{MapSC2} show different examples for mappings of carpets by $ \varPhi(x, y)=(x^2+y^2,\; x-y) $ and $ \varPhi(x, y)=(\sin x+y,\; \cos x) $ respectively. The colors displayed in these figures are produced in a similar way to those in Fig. \ref{SLikeCarpet}.


\subsection{\normalsize Mapping of gaskets}

\begin{figure}[]
	\centering
	\subfigure[Image of Sierpinski gasket]{\includegraphics[width = 2.2in]{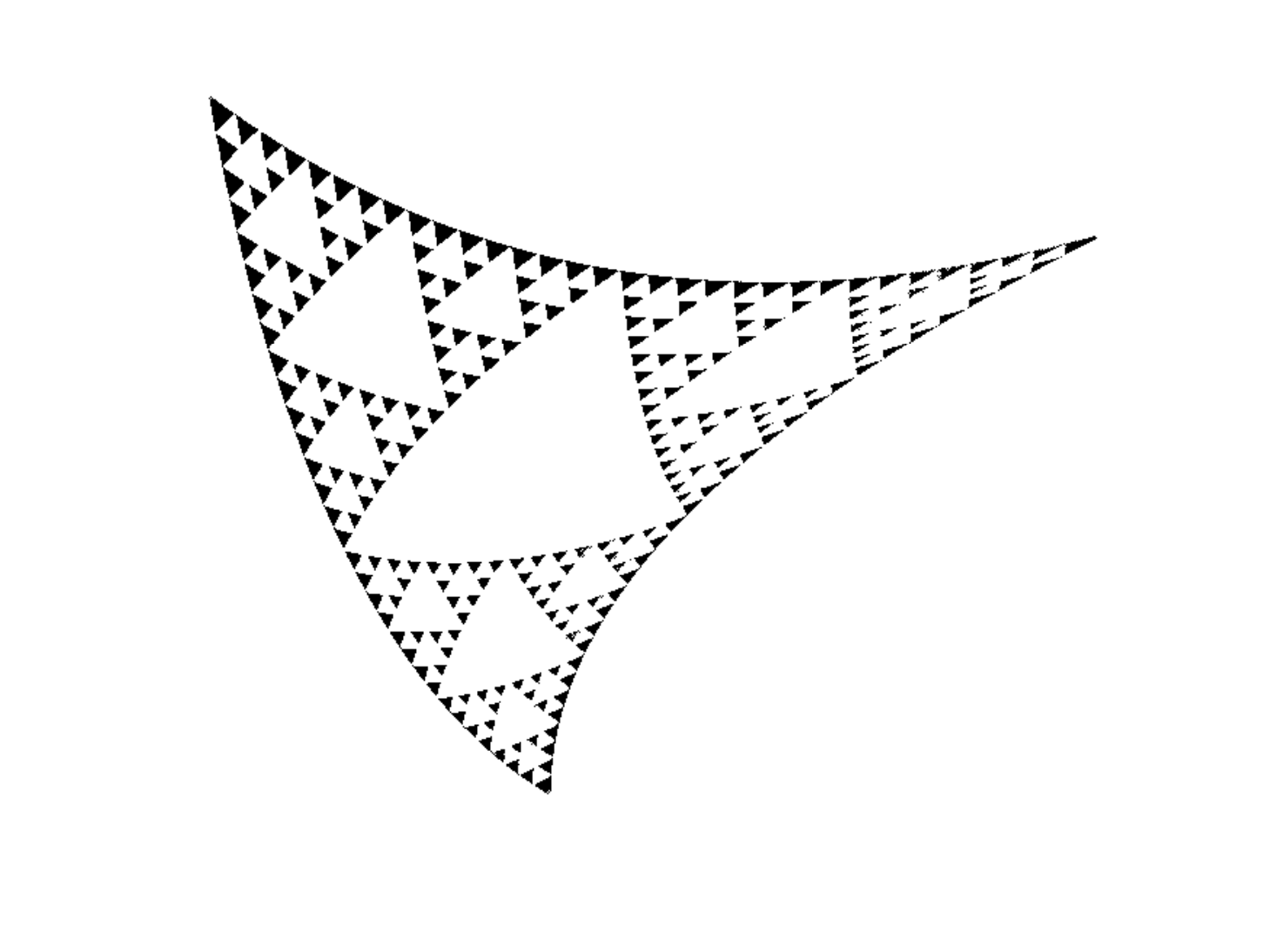}\label{}} 
	\subfigure[Image of gasket with $ \alpha(x)= \beta(x)= \gamma(x)= \cos x, \; a=2 $]{\includegraphics[width = 2.2in]{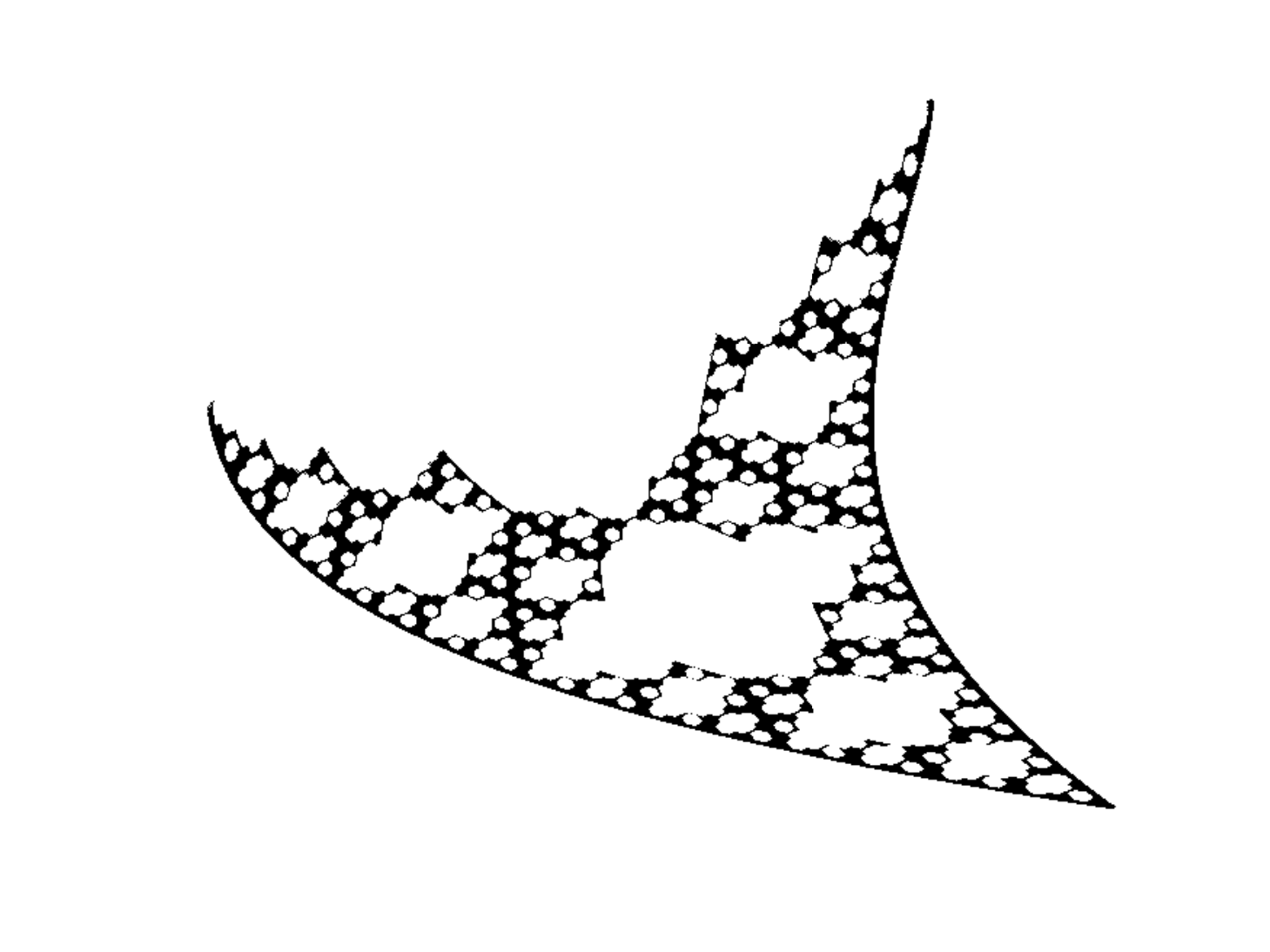}\label{}}
	\caption{Mappings of gaskets}
	\label{MapSG}   				
\end{figure}

Let us give examples of mappings of gaskets using formula (\ref{SGI}). Figure \ref{MapSG} (a) shows the mapped Sierpinski gasket by the map $ \varPhi(x, y)=(x^2-y, x+y^2) $, whereas Fig. \ref{MapSG} (b) represents the mapped gasket depicted in Fig. \ref{SLikeGasket} (b) by the map $ \varPhi(x, y)=(x+y^2, x-2y^{\frac{2}{3}}) $.

\section{\normalfont DYNAMICS}

Based on the fractal mapping iteration, we introduce dynamics in fractals. As in paper \cite{Akhmet}, discrete dynamics can be constructed for Sierpinski fractals by iterating the mappings introduced in the previous section. In other words, if we start, for example, with the Sierpinski carpet as an initial set, $ \mathcal{S}_0 $, and iterate the map $ \varPhi $ in FMI, we shall have the image sets
\[ \mathcal{S}_m=\varPhi^m(\mathcal{S}_0), \; m=1, 2, 3, ... \, . \]
Thus the discrete dynamics will consist of these fractal sets, $ \mathcal{S}_m $, as points of a trajectory. 

For continuous dynamics, the idea is to use the motion of a dynamical system with a fractal as an initial set. The motion of dynamical system is defined by $ \mathcal{A}_t \textbf{x}_0=\varphi(t, \textbf{x}_0) $, where $ \varphi $ is the solution of a two dimensional system of ordinary differential equations
\begin{equation} \label{ODE}
\textbf{x}'=g(t, \textbf{x}),
\end{equation}
with $ \varphi(0, \textbf{x}_0) = \textbf{x}_0 $. 

In the case of the Sierpinski carpet, we iteratively apply a motion $ \mathcal{A}_t $ to the scheme (\ref{SCI}) in the way
\begin{equation*}
\mathcal{A}_{-t}(\xi_n, \eta_n)= \psi_n \big( \mathcal{A}_{-t}(\xi, \eta) \big),
\end{equation*}
where $ \mathcal{A}_t(x, y)=(A_t x, B_t y) $. Through this procedure, we construct dynamics of sets  $A_t \mathcal{F},$ where the Sierpinski carpet $ \mathcal{F} $ is the initial value. Thus, the differential equations are involved in fractals such that the latter become points of the solution trajectory. If the map $ A_t $ is bi-Lipschitzian (this is true, for instance, if the function $ g $ in (\ref{ODE}) is Lipschitzian) then the  set $ A_t \mathcal{F} $ for each fixed $ t $ is a fractal.

Let us now consider the Van Der Pol equation
\begin{equation} \label{VanDPEq}
 u''+\mu(u^2-1)u'+u=0,
\end{equation}
where $ \mu $ is a real constant known as the damping parameter. Using the variables $ x=u $ and $ y=u' $, one can show that Eq. (\ref{VanDPEq}) is equivalent to the autonomous system
\begin{equation} \label{VDP System}
\begin{split}
&x' = y, \\
&y'= \mu(1-x^2) y-x.
\end{split}	
\end{equation}

\begin{figure}[]
	\centering
	\includegraphics[width=0.5\linewidth]{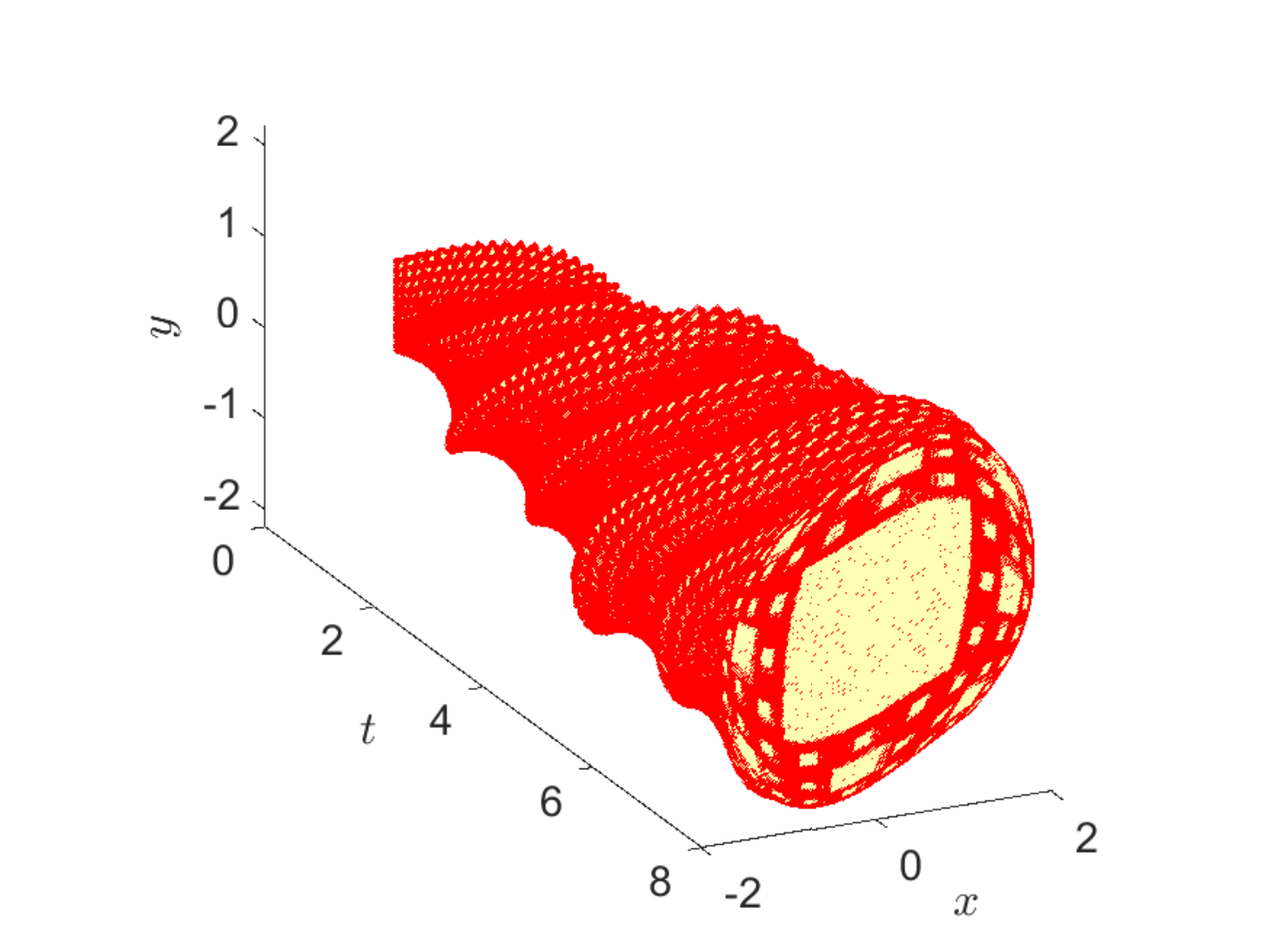}
	\caption{Van Der Pol dynamics of Sierpinski carpet}
	\label{VDPDynSC1}	
\end{figure}

\begin{figure}[H]
	\centering	
	\begin{minipage}{0.9\textwidth}
		\centering	
		\subfigure[$ t=1 $]{\includegraphics[width = 1.8in]{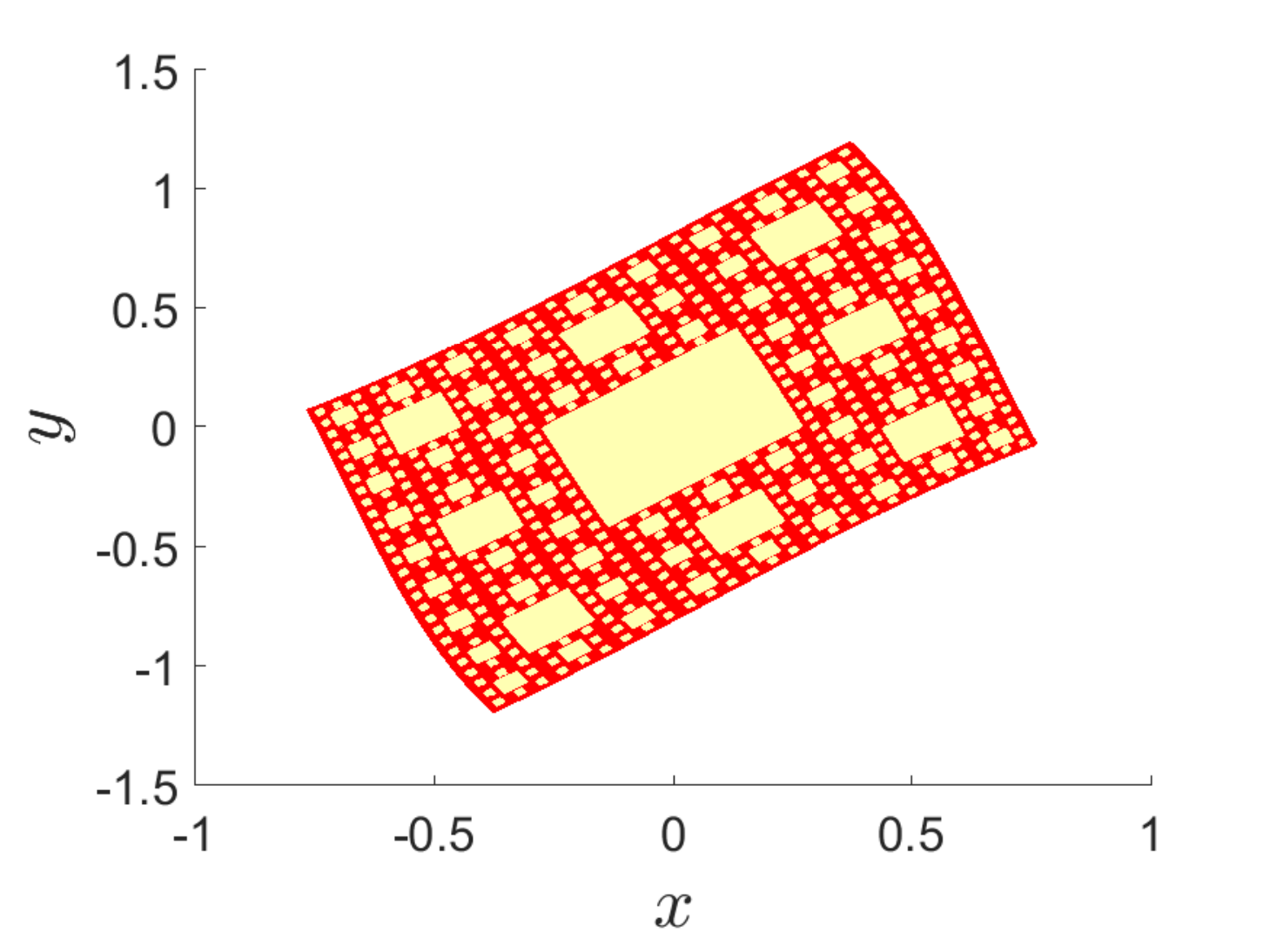}} \hspace{1.5cm}
		\subfigure[$ t=3 $]{\includegraphics[width = 1.8in]{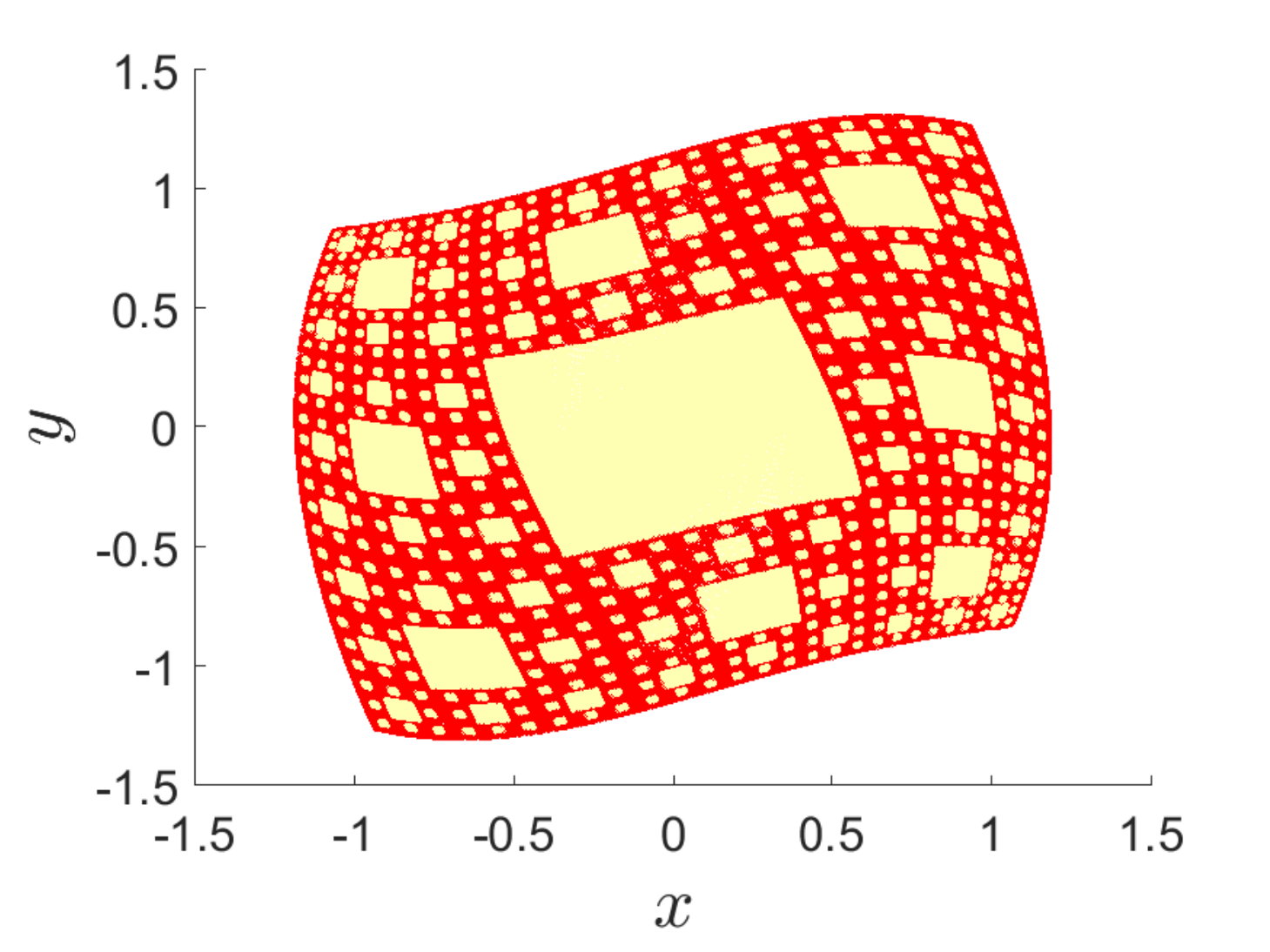}}
	\end{minipage} 
	\begin{minipage}{0.9\textwidth}	
		\centering		
		\subfigure[$ t=5 $]{\includegraphics[width = 1.8in]{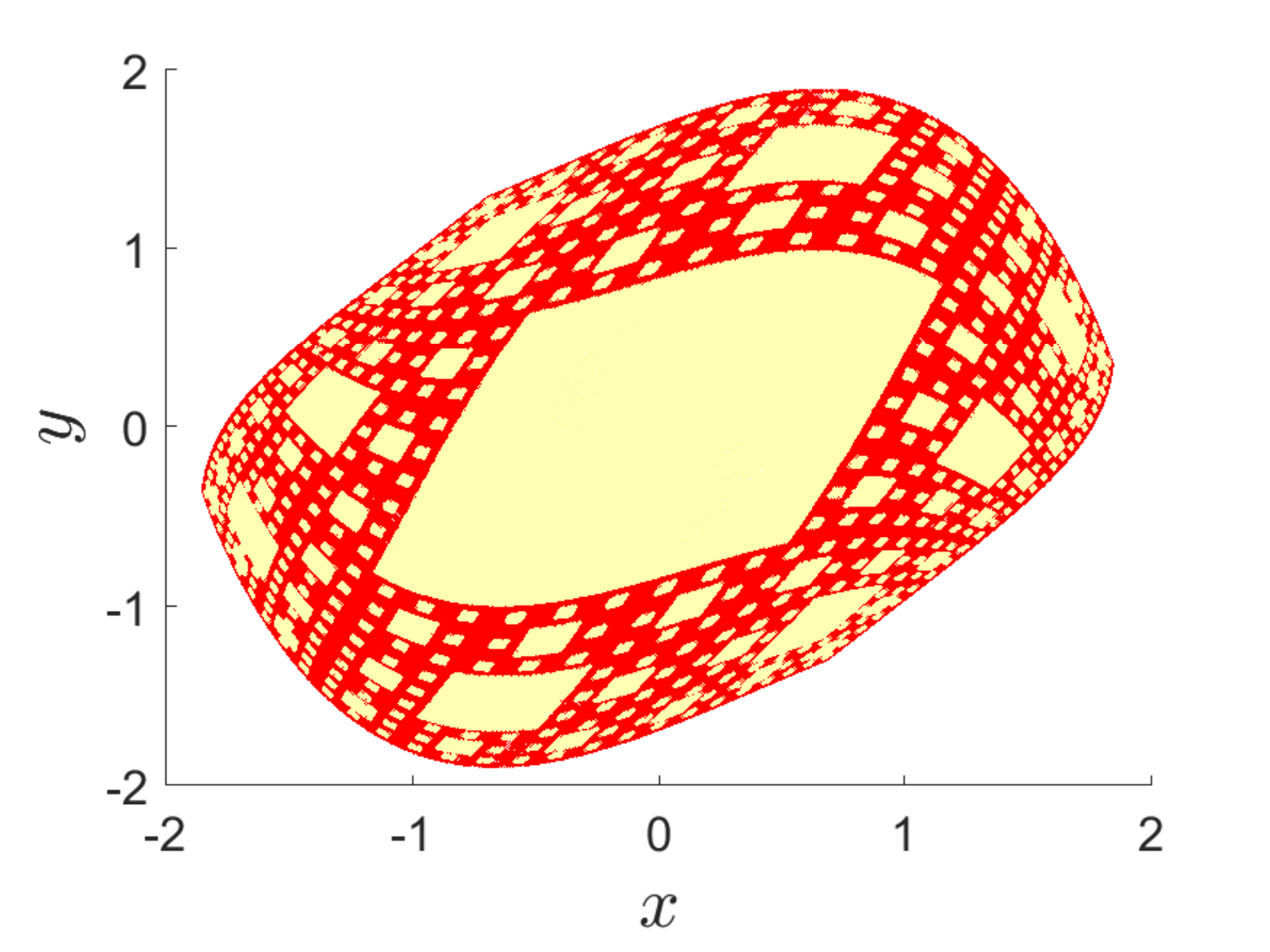}}\hspace{1.5cm}
		\subfigure[$ t=7 $]{\includegraphics[width = 1.8in]{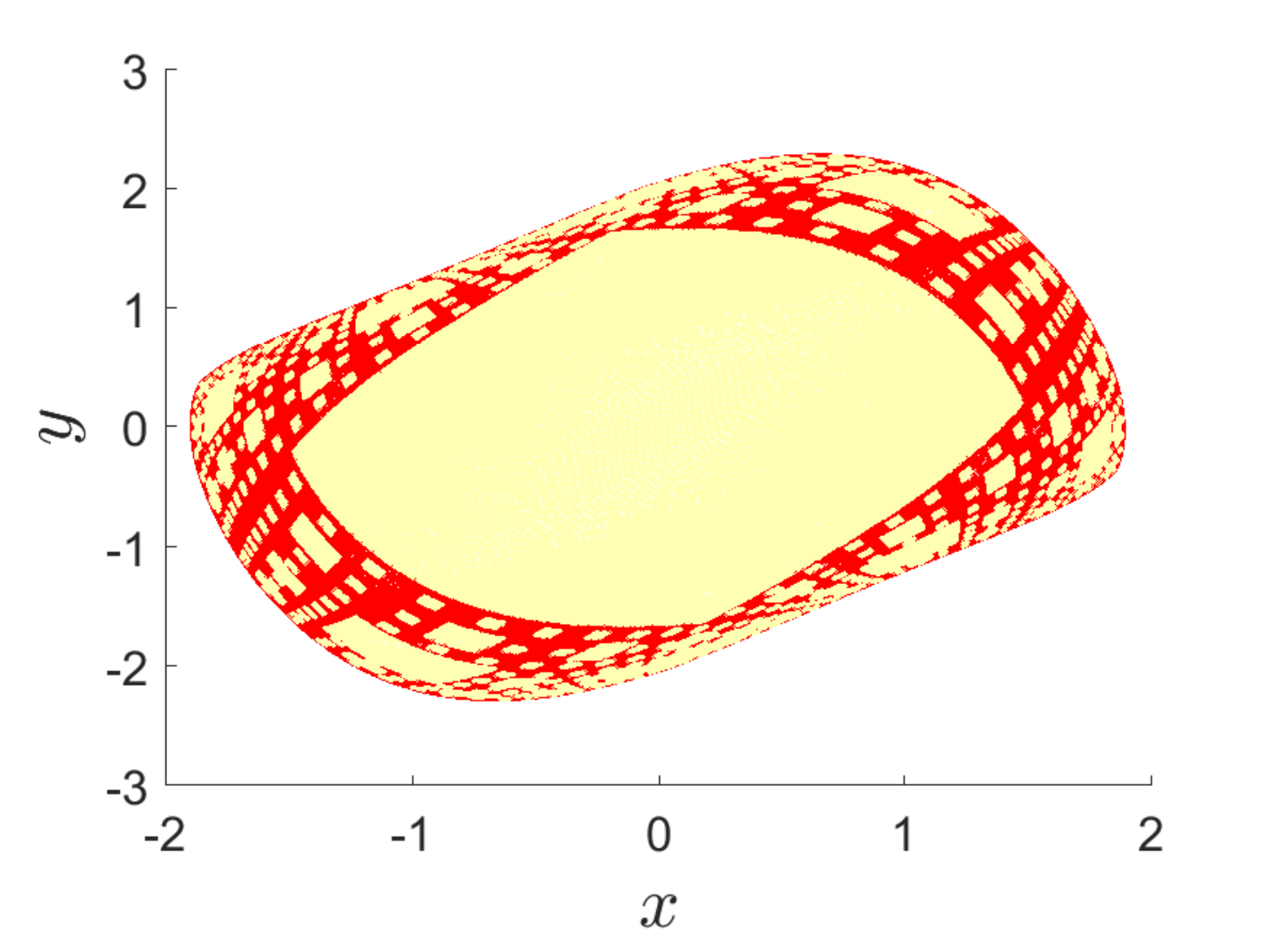}}
	\end{minipage} 
	\caption{Trajectory sections of the Van Der Pol dynamics in Sierpinski carpet}
	\label{VDPDynSCSections}   				
\end{figure}

Let us denote by $ \big(x(t, x_0), y(t, y_0)\big) $ the solution of (\ref{VDP System}) with $ x(0, x_0)= x_0, \; y(0, y_0)= y_0 $. System (\ref{VDP System}) can be numerically solved to construct a dynamical system with the motion $ \mathcal{A}_t(x_0, y_0)=(A_t x_0, B_t y_0) $ where $ A_t x_0=x(t, x_0) $ and $ B_t y_0=y(t, y_0) $. We apply this dynamics for an approximation of the Sierpinski Carpet as an initial set. The trajectory of the Van Der Pol dynamics with $ \mu=0.5 $ and $ 0 \leq t \leq 8 $ is shown in Fig. \ref{VDPDynSC1}. Figure \ref{VDPDynSCSections} exhibits the sections of the trajectory at the moments $ t=1, t=3, t=5 $, and $ t=7 $.

\begin{figure}[]
	\centering
	\subfigure[$ t=1 $]{\includegraphics[width = 2.0in]{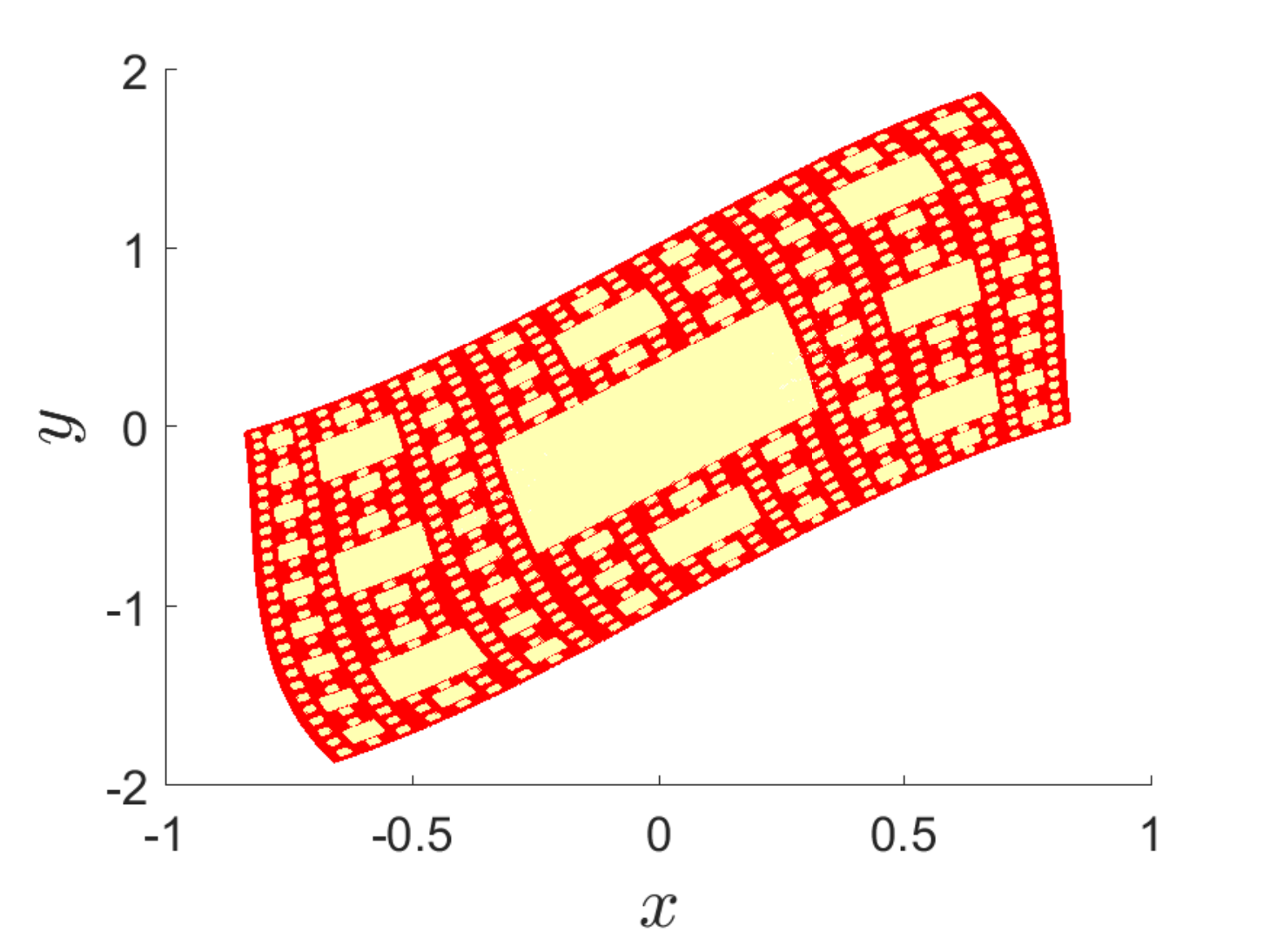}\label{}} \hspace{1.5cm}
	\subfigure[$ t=3 $]{\includegraphics[width = 2.0in]{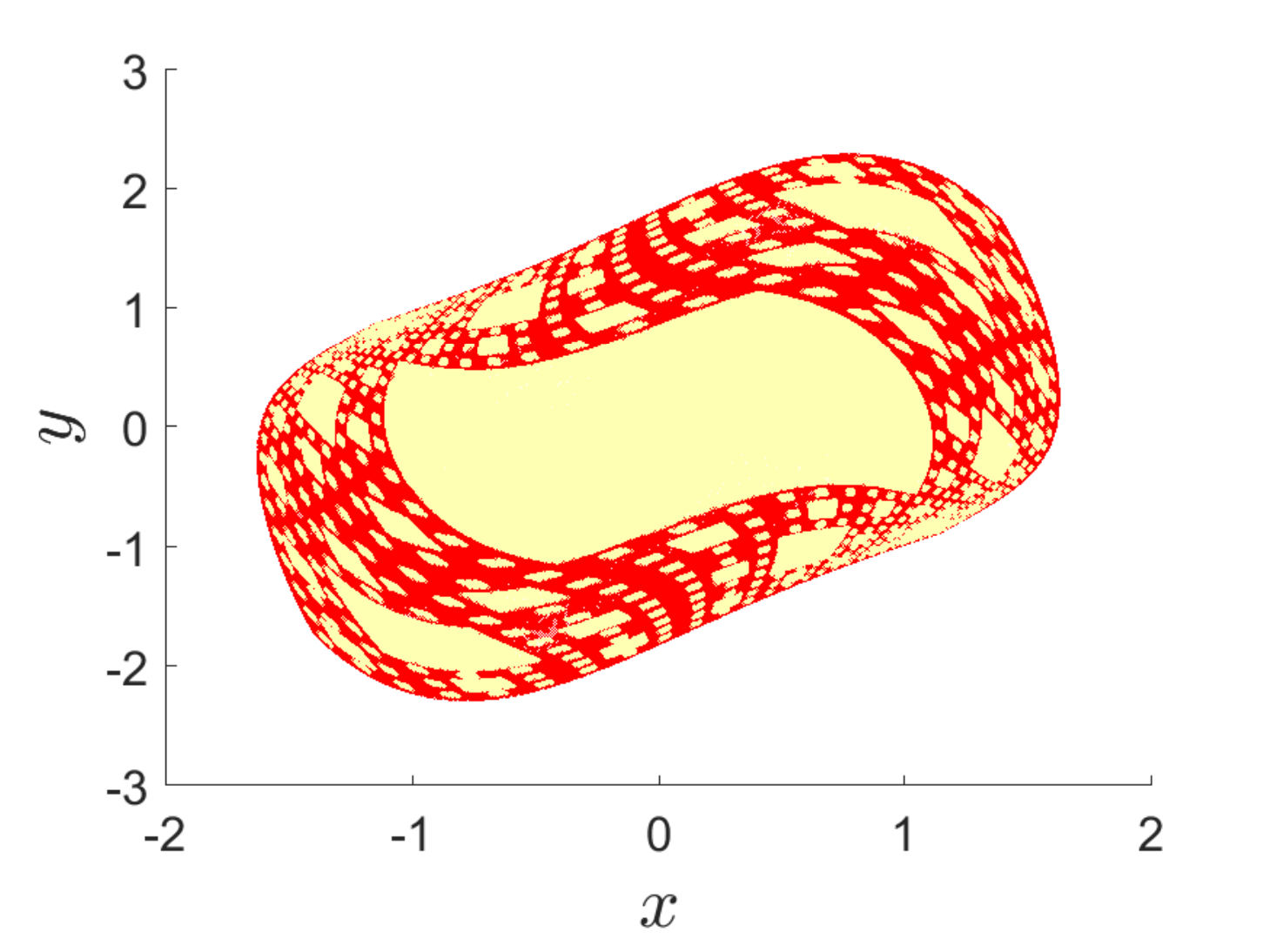}\label{}}
	\caption{Trajectory sections of the Van Der Pol dynamics with $ \mu=1.3 $}
	\label{VDPDynSCSec-a}   				
\end{figure}

In Fig. \ref{VDPDynSCSec-a} we again show two sections of the trajectory of the Van Der Pol dynamics for the approximation of the Sierpinski Carpet but with $ \mu=1.3 $. Comparing these sections with their time counterpart in Fig. \ref{VDPDynSCSections}, we can observe the dissimilarity in the deformation rate in the structure of the Sierpinski carpet. This is attributed to that the value of the damping parameter reflects the degree of nonlinearity of the Van Der Pol equation.

For dynamics in Sierpinski gasket, let us consider the Duffing equation
\[ u''+\delta u'+\beta u+\alpha u^3=\gamma \cos \omega t, \]
where $ \delta, \beta, \alpha, \gamma $, and $ \omega $ are real parameters. The equation is equivalent to the non-autonomous system
\begin{equation*} \label{Duffings System}
\begin{split}
&x' = y, \\
&y'= -\delta y-\beta x-\alpha x^3+\gamma \cos \omega t.
\end{split}	
\end{equation*}

\begin{figure}[]
	\centering
	\includegraphics[width=0.5\linewidth]{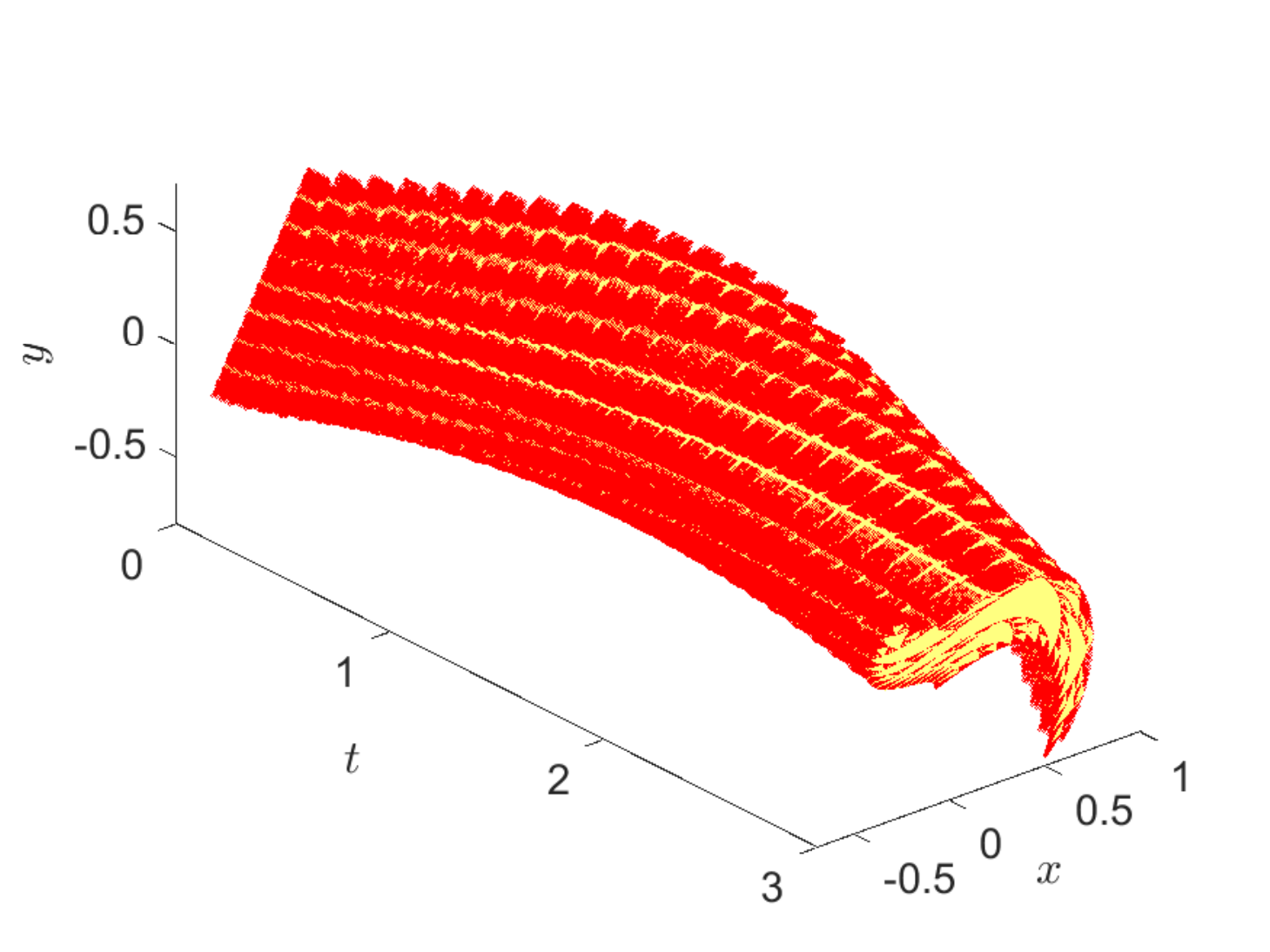}
	\caption{Trajectory of the Duffing dynamics in Sierpinski gasket}
	\label{DuffDynSG}	
\end{figure}

\begin{figure}[]
	\centering	
	\begin{minipage}{0.9\textwidth}
		\centering	
		\subfigure[$ t=0.8 $]{\includegraphics[width = 1.8in]{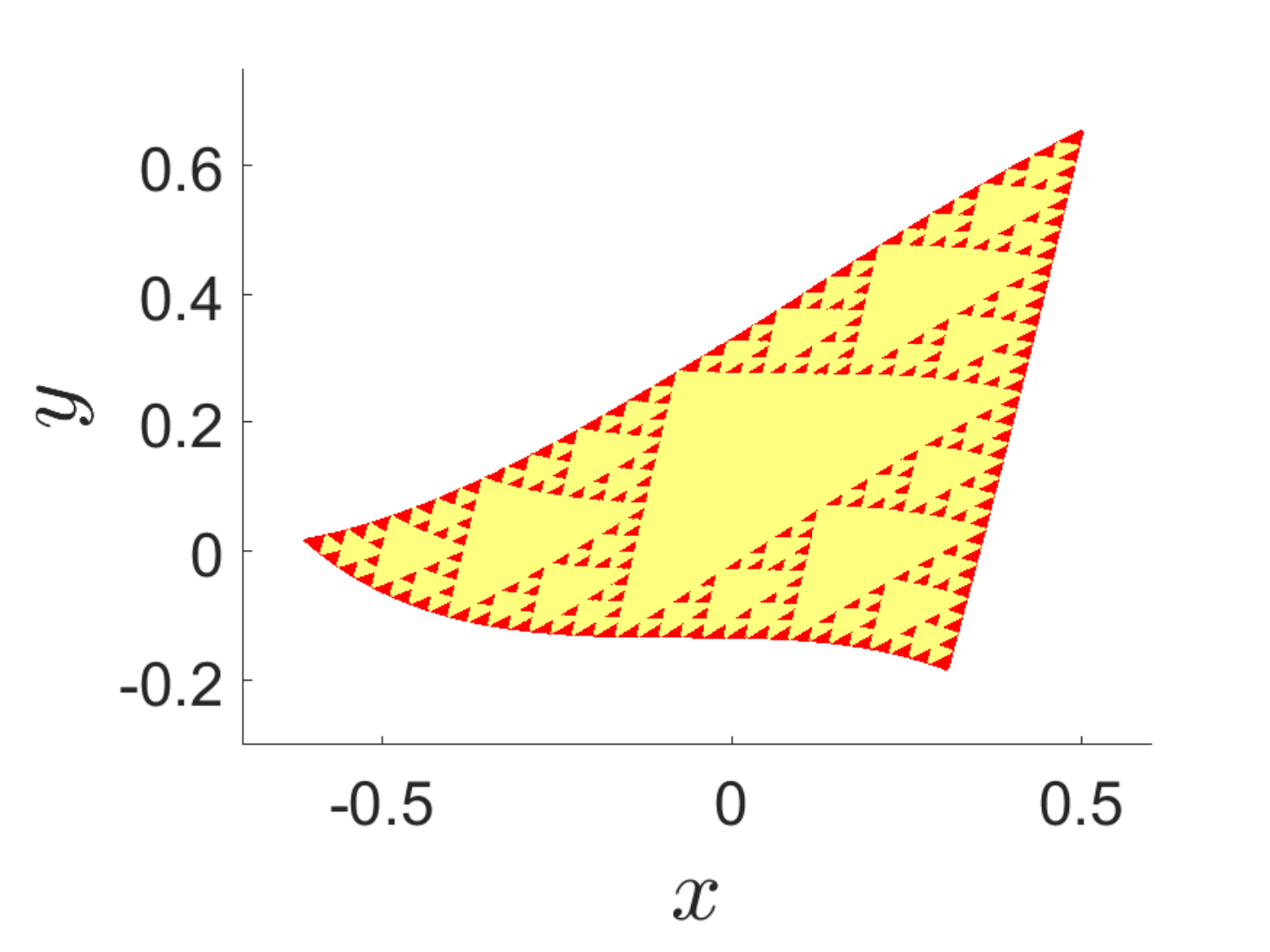}} \hspace{1.5cm}
		\subfigure[$ t=1.4 $]{\includegraphics[width = 1.8in]{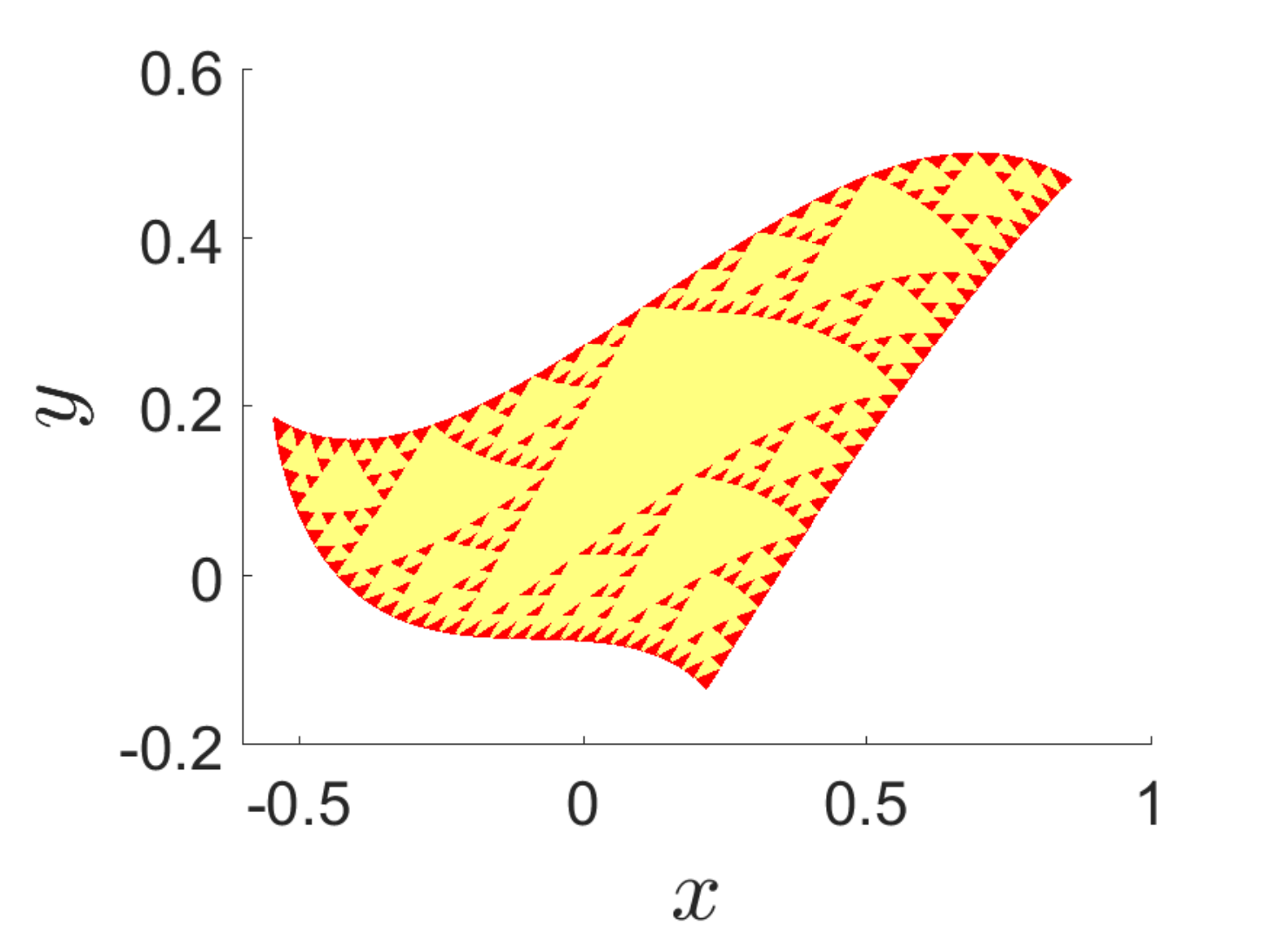}}
	\end{minipage} 
	\begin{minipage}{0.9\textwidth}	
		\centering		
		\subfigure[$ t=2.0 $]{\includegraphics[width = 1.8in]{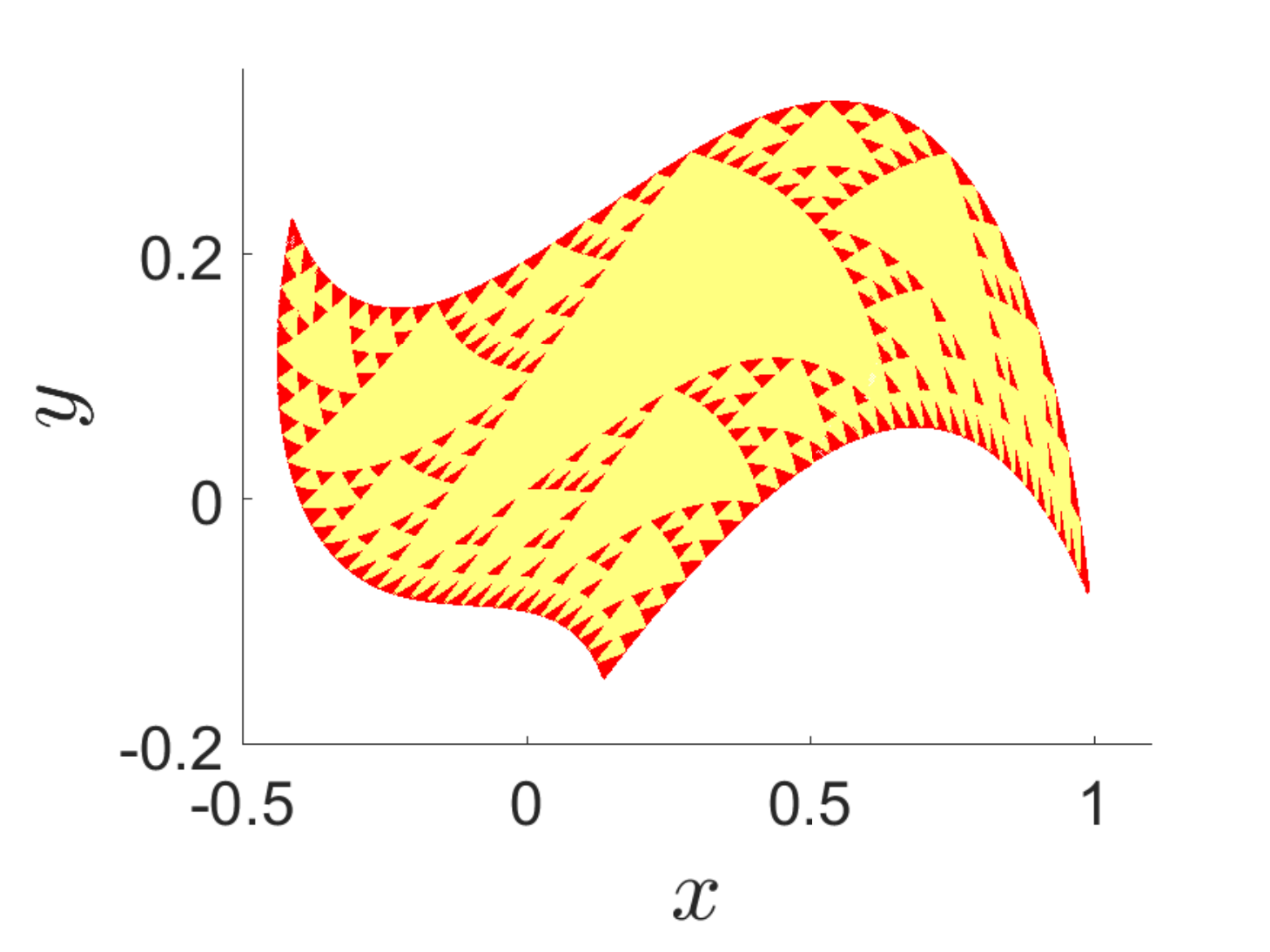}}\hspace{1.5cm}
		\subfigure[$ t=2.6 $]{\includegraphics[width = 1.8in]{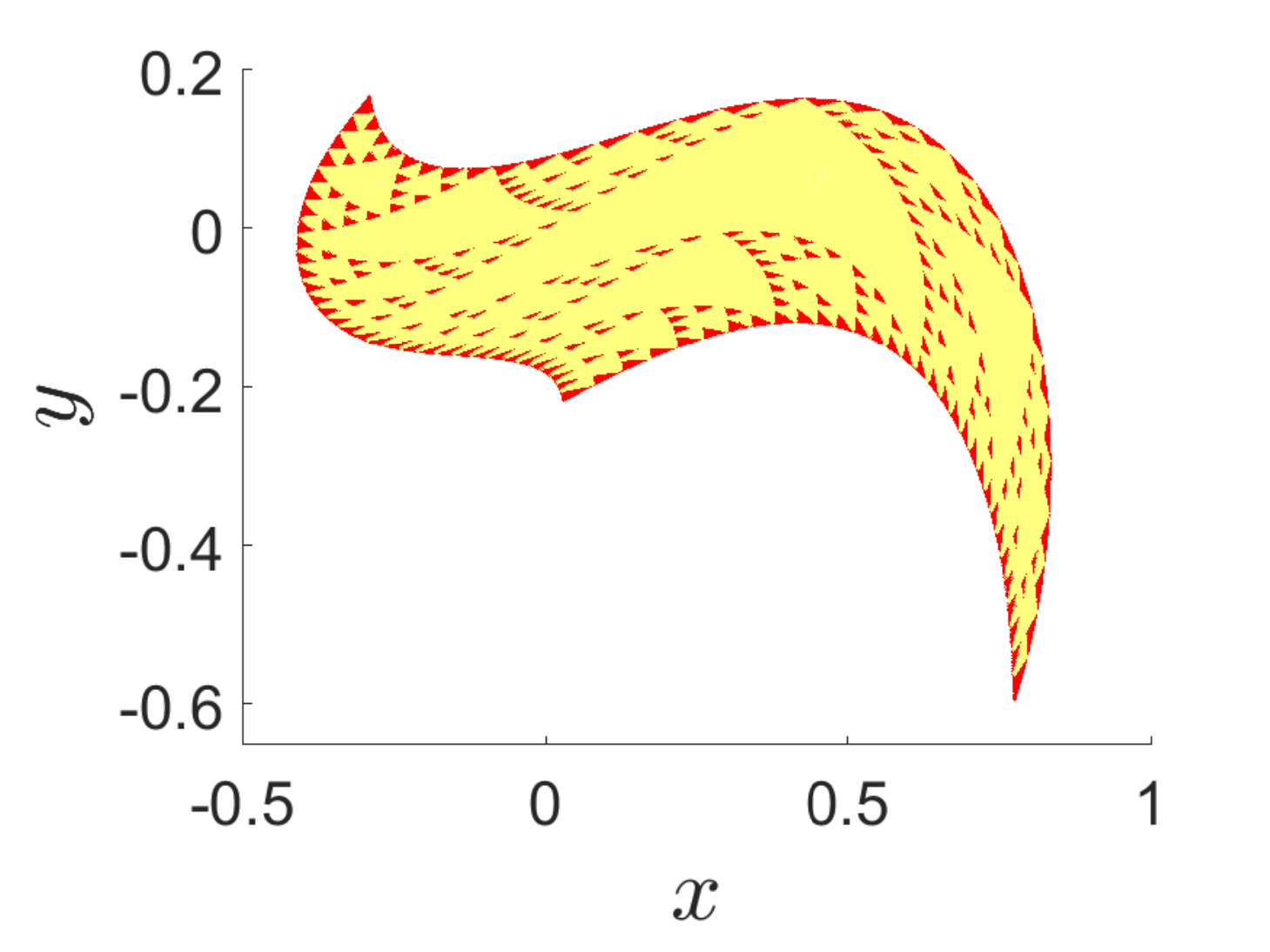}}
	\end{minipage} 
	\caption{Sections of the trajectory}
	\label{DuffDynSGSections}   				
\end{figure}

In a similar way to the mapping of gasket, we apply the dynamical system associated with the Duffing equation to an approximation of the Sierpinski gasket. The fractal trajectory for $ 0 \leq t \leq 3 $ is shown in Fig. \ref{DuffDynSG}, whereas Fig. \ref{DuffDynSGSections} displays the sections of the trajectory at the specific times $ t=0.8, \, t=1.4, \, t=2.0 $, and $ t=2.6 $. The values $ \delta=0.08, \, \beta=0, \, \alpha=1, \, \gamma=0.2 $ and $ \omega=1 $ are used in the simulation.

\section{\normalfont CONCLUSION}

In the present paper we consider the dynamics for Sierpinski fractals. By the dynamics we mean:
\begin{enumerate}[label=(\roman*)]
	\item Dynamics of iterations used for generating fractals;
	\item Dynamics of iterations that map fractals;
	\item Construction of discrete trajectories of fractals by iterating the map;
	\item Construction of continuous fractal trajectories through differential equations and iterations.
\end{enumerate}

The main purpose of this study is to develop and extend the idea of FMI which have been introduced in our previous paper \cite{Akhmet}. In that paper the FMI was applied only for the Julia and Mandelbrot sets. Here we focus on Sierpinski fractals, how to generate them through iteration schemes, and how to formulate FMI for these schemes in order to construct discrete and continuous dynamics in the fractals.

Owing to the important roles of Sierpinski fractals in several applications like weighted networks, trapping problems, antenna engineering, city planning, and urban growth \cite{Dong,Zhang,Kozak,Kansal,Triantakonstantis}, we expect that the results of the present study will be helpful in the fields of applications. One of the crucial applications of fractals  involves optimization theory. Fractal geometry is used to solve some classes of optimization problems such as supply chain management and hierarchical design \cite{Kirkhope,Liu}. In paper \cite{Kirkhope}, for instance, the properties of a particular hierarchical structure is established. The authors constructed the relationship between the Hausdorff dimension of the optimal structure and loading for which the structure is optimized. The Hausdorff dimension is calculated through considering the self-similarity of the structure at different hierarchical levels. The self-similar fractals, like the Sierpinski gasket, are considered as effective tools for studying the hierarchical structures \cite{Kozak,Riera}. Thus, finding a way to map this type of structures allows to create a new hierarchical structure with the same Hausdorff dimension but different mechanical properties if one consider bi-Lipschitz maps.

Further applications can be considered by taking into account the relationship between the fractal theory of motion and quantum mechanics. In the scale relativity theory \cite{Nottale1,Nottale2}, fractals are considered as a geometric framework of atomic scale motions such that the quantum behavior can be viewed as particles moving on fractal trajectories. One can suppose that by composing the scale relativity theory with dynamics of fractals developed in this paper, we will be able to understand better the fractal nature of the world.  A possible connection between fractal mappings and quantum mechanics through the scale relativity theory can provide important applications for the the former in various fields such as biology, cosmology, and fractal geodesics (see \cite{Nottale2} and the relevant references therein). 

\section*{\normalfont ACKNOWLEDGMENTS}

The third author is supported by a scholarship from the Ministry of Education, Libya.

\end{document}